\documentclass[preprint,12pt]{elsarticle}

\usepackage{amssymb}
\usepackage{amsmath}
\usepackage{amsfonts}
\usepackage{amsthm}
\usepackage{algorithm}
\usepackage{algpseudocode}
\usepackage[shortlabels]{enumitem}
\usepackage{url}
\usepackage{bm}
\usepackage{mathrsfs}
\usepackage{caption}
\usepackage{subcaption}
\usepackage{booktabs}
\usepackage[group-separator={,}]{siunitx}

\usepackage{graphicx}
\usepackage{caption}
\usepackage{xcolor}
\usepackage{csquotes}
\usepackage{pgfplots} 
\pgfplotsset{compat=1.15} 
\usetikzlibrary{calc}

\definecolor{codeorange}{rgb}{0.9,0.65,0}
\definecolor{codegreen}{rgb}{0,0.6,0}
\definecolor{codegray}{rgb}{0.5,0.5,0.5}
\definecolor{codepurple}{rgb}{0.58,0,0.82}
\definecolor{carnelian}{rgb}{0.7, 0.11, 0.11}
\definecolor{backcolour}{rgb}{0.95,0.95,0.92}

\algnewcommand\algorithmicforeach{\textbf{for each}}
\algdef{S}[FOR]{ForEach}[1]{\algorithmicforeach\ #1\ \algorithmicdo}

\newcommand{\dx}[1]{\mathrm{d}{#1}}
\newcommand{\epsrel}{\epsilon_{\text{rel}}}
\newcommand{\nn}{n_{\text{nodes}}}
\newcommand{\nf}{n_{\text{facets}}}
\newcommand{\Jsparse}{J^{\text{sparse}}}
\newcommand{\Jcav}{J^{\text{cav}}}
\newcommand{\xtimes}[1]{{#1}$\times$}

\newlist{romanitemize}{enumerate}{1}

\newtheorem{example}{Example}

\setlist[romanitemize, 1]
{label=\alph{romanitemizei}., 
leftmargin=30pt,
rightmargin=10pt
}

\journal{CMAME}

\begin{document}

\begin{frontmatter}



\title{Hierarchical Block Low-rank Approximation \newline of  Cavity Radiation}

\author[eth,akselos]{Ivan Baburin}
\ead{ivan.baburin@inf.ethz.ch}

\author{Jonas Ballani}
\ead{jonas.ballani@gmail.com}

\author[akselos]{John W. Peterson}
\ead{john.peterson@akselos.com}

\author[akselos]{David Knezevic\corref{cor1}}
\ead{david.knezevic@akselos.com}

\cortext[cor1]{Corresponding author}

\affiliation[eth]{organization={ETH},
city={Z\unexpanded{\"u}rich},
country={Switzerland}
}

\affiliation[akselos]{organization={Akselos SA},
city={Lausanne},
country={Switzerland}
}

\begin{abstract}
In this paper we examine the use of low-rank approximations for the handling of cavity radiation boundary conditions in a transient heat transfer problem. The finite element discretization that arises from cavity radiation involves far-field degree of freedom coupling and matrices characterized by a mixed dense/sparse structure, which pose difficulties for the efficiency and scalability of solvers. Here we consider a special treatment of the cavity radiation discretization using a block low-rank approximation combined with hierarchical matrices. We provide an overview of the methodology and discuss techniques that can be used to improve efficiency within the framework of hierarchical matrices, including use of the adaptive cross approximation (ACA) method. We provide a number of numerical results that demonstrate the accuracy and efficiency of the approach in practical problems, and demonstrate significant speed-up and memory reduction compared to the more conventional ``dense matrix'' approach.
\end{abstract}




\end{frontmatter}


\section{Introduction\label{sec:intro}}

Blackbody radiation is a well-studied phenomenon in classical mechanics. Though challenging to model mathematically (c.f.\ the ``paradox'' known as the ``ultraviolet catastrophe'' \cite{McQuarrie_1997}), the study of blackbody radiation at the beginning of 20th century gave birth to Planck's law and the concept of quantization of energy, and led to further fundamental discoveries that laid the foundations for the field of quantum theory. In the modern era of scientific computation, there are still many challenges connected to blackbody radiation, in particular the computation of radiation-based heat transfer. In this article we consider a more general setting of ``graybody'' radiation, in which the emissivity coefficient may be specified on the cavity surface. This is relevant in many industrial and scientific settings, such as heat transfer in reactors \cite{Chinoy_1991,vonZedtwitz_2005,ZGraggen_2008}, or radiative heat transfer on celestial surfaces \cite{Barman_2001,Delbo_2005,Hayne_2021,Potter2022}.

The primary challenge of cavity radiation from the computational point of view is that the interactions are ``non-local'': in general every surface facet interacts with every other surface facet, where the interactions are quantified via the so-called ``reflection matrix,'' which in general is dense and of size $n_{\rm facets} \times n_{\rm facets}$. This non-locality results in a loss of sparsity in the matrices arising from finite element discretizations, which in turn leads to very high computational cost when $n_{\rm facets}$ becomes large (e.g.\ $>\num{10000}$). The most straightforward approach --- explicitly allocating and computing the dense matrices associated with cavity radiation --- is prohibitively expensive in terms of both memory utilization and computation time, especially when applied on large or highly detailed models.

We propose an alternative method of computing cavity radiation which avoids dense matrix discretizations by exploiting the physical nature of radiation, namely its fast (inverse-square) decay relative to the distance between emitting surfaces. Due to this property, many individual entries in the view factor matrix $F$ and reflection matrix $C$ (for the precise definitions we refer to Section \ref{sec:problem_setting}) underlying the computation of the cavity radiation flux interact only weakly, and thus, if grouped together, can be approximated reliably using low-rank blocks. Since the grouping or clustering of facets is strongly linked to the underlying geometry, it is natural to store the $F$ and $C$ matrices hierarchically, identifying large, weakly interacting matrix blocks through a geometric clustering criterion. In the analogous context of boundary element methods, hierarchical matrices have been successfully applied to large-scale problems, see \cite{Hackbusch2015, Sauter_2011} for comprehensive introductions.


While we were working on the implementation, a very similar approach was presented by Potter et al.~\cite{Potter2022} in the context of thermal irradiance on planetary surfaces, where the radiosity matrix shares many structural properties with the view factor matrix $F$ from heat transfer. In their paper, Potter et al.\ describe an efficient construction of a hierarchical matrix with low-rank approximation of individual blocks, and provide an application of their approach to the computation of irradiance on a lunar crater. The numerical results demonstrate the potential benefits of the hierarchical approach when applied to practical problems.

The novelty of our approach with regard to Potter et al.\ is twofold. First, due to the nature of thermal radiation, we will perform computations not only with the hierarchical view factor matrix $F$, but also with the inverse of the reflection matrix $C$. Due to the nonlinear (temperature-dependent materials, radiation flux terms) and time-dependent character of the heat equation, we have to perform many applications of $C^{-1}$, but, because $C$ is itself \emph{independent} of the current temperature solution, we are able to employ a \textit{block-LU decomposition} of $C$. While there is a non-trivial up-front computational cost involved in the construction of the block-LU decomposition, this cost is amortized over the many time steps and nonlinear iterations which are performed during the course of the simulation. We note that the low-rank structure of the reflection matrix is preserved in its block-LU decomposition, and this fact is exploited to allow the action of $C^{-1}$ to be computed quickly and with low memory requirements.

Second, we use the \textit{adaptive cross approximation} (ACA) \cite{Bebendorf2000} to optimize the construction of the view factor matrix, since we found that it was significantly faster than the singular value decomposition when constructing a low-rank approximation of individual blocks. Due to its strong reliability guarantees, we have opted to use ACA with \textit{full pivoting} (running linearly in the number of entries). Further computational gains can be expected when using ACA with \textit{partial pivoting} (running linearly in the number of facets), at the cost of potentially less reliable low-rank approximations due to the use of partial information.
A comprehensive discussion of different cross approximation techniques can be found in \cite{Borm_2006}.

Note also that in the current work we do not consider {\em occlusion detection} of surface facets, and hence in Section~\ref{sec:numerical_experiments} we conduct numerical experiments in which occlusion detection is not required. Including occlusion detection is a natural extension of the implementation, e.g.\ via a ray-tracing approach \cite{Franklin_thesis}. We emphasize that incorporating occlusion detection would not change the block-low rank approach that we present here in any way, since it corresponds to ``zeroing out'' entries in the view factor matrix that correspond to occluded facets, which is equivalent to modifying only the ``input'' to the block low-rank algorithm.




\section{Problem Setting\label{sec:problem_setting}}

\subsection{Transient heat equation}
We consider transient heat transfer in a domain $\Omega \subset \mathbb{R}^3$ with radiation boundary $\Gamma \subseteq \partial \Omega$ governed by:
\begin{alignat}{2}
    \label{eq:heat}
    \rho c_p \frac{\partial T}{\partial t}   & = \nabla \cdot (k \nabla T) + f && \quad \in \Omega 
    \\
    \label{eq:radiation_bc}
    k \nabla T \cdot \vec{n} & = q(T) && \quad \in \Gamma 
\end{alignat}

\noindent where $T$ is the unknown temperature, $q(T)$ is the radiation heat flux, $f$ is the internal heat source, $k$ is the thermal conductivity, $\rho$ is the material density, and $c_p$ is the specific heat. In addition to the radiation boundary condition \eqref{eq:radiation_bc} shown above, there may in general be other types of boundary conditions (Dirichlet, Neumann) present, though their details are not crucial to our discussion. Furthermore, in order to fully specify the problem, an initial condition for the temperature which applies on all of $\Omega$ at time $t=0$ must also be provided. 

The variational statement associated with \eqref{eq:heat} is then: find $T \in \mathcal{H}^1(\Omega)$ such that

\begin{equation}\label{eq:SVF}
    \int_{\Omega} \rho c_p \frac{\partial T}{\partial t} v \, \dx{V} + \int_{\Omega} k \nabla T \cdot \nabla v \, \dx{V} - \int_{\Gamma} q(T)v \, \dx{S} = \int_{\Omega} f v \, \dx{V} 
\end{equation}

\noindent holds for all test functions $v$ from the Sobolev space $\mathcal{H}^1(\Omega)$.
Equation \eqref{eq:SVF} is in general a nonlinear equation for $T$, both because the flux $q(T)$ depends on $T$ to the fourth power (see \eqref{eq:emission}), but also because the thermal conductivity and specific heat may depend on the temperature. The nonlinearity of \eqref{eq:SVF} plays a crucial role in its numerical solution, since we need to iterate (via Newton's method) within each time step to find the converged solution at the current time, and each iteration requires the formation and solution of a large linear system of equations involving the Jacobian matrix, as we will further elaborate in Section~\ref{sec:discretization_details}.

\subsection{Radiation flux\label{sec:rad_flux}}
In this section, we briefly recall some of the fundamentals of thermal radiation which are relevant to the current work, as presented in \cite{Holman1986}. The thermal energy flux (amount of energy emitted per unit time, per unit area) from a so-called ``graybody'' thermal radiator is given, according to the Stefan--Boltzmann law of radiation, by

\begin{equation}
    \label{eq:emission}
    q = \epsilon \sigma T^4 
\end{equation}

\noindent
where $\sigma \approx 5.669 \times 10^{-9} \; \frac{W}{m^2 K^4}$ is the Stefan--Boltzmann constant, $0 < \epsilon \leq 1$ is the surface emissivity (dimensionless), and $T$ is the absolute surface temperature of the body, measured in Kelvin. As $\epsilon \rightarrow 1$, we recover the case of a ``blackbody'' radiator. In this formulation, we also assume that the emissivity is independent of the wavelength of the radiation, that is, the surfaces are so-called ``ideal'' graybodies \cite{Chapman_1984}. The SI units of the energy flux $q$ are thus $\frac{W}{m^2}$. 

In the case of two radiating graybody infinite flat planes with temperatures $T_1$ and $T_2$ and surface emissivities $\epsilon_1$ and $\epsilon_2$, respectively, it can be shown \cite{Chapman_1984} that the net radiation energy flux between them is proportional to the difference between their temperatures raised to the fourth power:

\begin{equation}
    q = \frac{ \sigma \epsilon_1 \epsilon_2 (T_1^4 - T_2^4) }{1 - (1-\epsilon_1)(1-\epsilon_2)}
\end{equation}
\noindent The case of infinite parallel planes is particularly simple because one can assume that all the energy emitted by plane 1 is absorbed by plane 2, and vice-versa, which then admits a closed-form solution based on Kirchoff's law. 

In real applications, one must of course consider finite, non-planar surfaces with non-uniform temperature fields. The standard approach taken in the general case is thus to discretize the surfaces comprising the radiation cavity into approximately planar finite regions or ``facets,'' and then consider the exchange of energy between individual pairs of facets separately, the details of which are discussed in \S\ref{sec:rad_flux_discrete}.

\section{Discretization details\label{sec:discretization_details}}

\subsection{Finite element discretization\label{sec:fe_discretization}}

First we discretize the three-dimensional domain $\Omega$ into a mesh $\mathcal{T}_h$ with $\nn$ nodes comprising hexahedral or tetrahedral elements with circumspherical diameter $\mathcal{O}(h)$. Our approach does not explicitly depend on the geometric type of the elements used, and should in theory work for all standard element types as well as ``hybrid'' meshes consisting of a mixture of geometric types. Then, the faces of the elements in $\mathcal{T}_h$ which lie on the radiation boundary $\Gamma$ can be considered as a lower-dimensional ``manifold'' mesh $\mathcal{M}_h$ consisting of $n_{\rm facets} := |\mathcal{M}_h|$ facets, which will all be allowed to emit and reflect radiation.

We next introduce the standard continuous Lagrange nodal basis $\{\varphi_{N}\}_{N=1}^{\nn}$ associated to $\mathcal{T}_h$ and state the semi-discrete Galerkin finite element approximation of \eqref{eq:SVF}, find $T_h$ such that:
\begin{equation}
    \label{eq:SVF_Galerkin}
      \int_{\Omega} \rho c_p \frac{\partial T_h}{\partial t} \varphi_{N} \, \dx{V} + \int_{\Omega} k \nabla T_h \cdot \nabla \varphi_{N} \, \dx{V} - \int_{\Gamma} q(T_h) \varphi_{N} \, \dx{S} = \int_{\Omega} f \varphi_{N} \, \dx{V} 
\end{equation}
holds for all $N=1,\ldots,\nn$, where
\begin{align}
    \label{eq:Th}
    T_h(x) := \sum_{M=1}^{\nn} T_{M} \varphi_{M}(x)
\end{align}
and $T_{M}$ are the unknown coefficients. The argument $x$ in \eqref{eq:Th} indicates the spatial dependence of the basis functions (and lack thereof for the coefficients). The temperature on node $M$ is consequently given by $T_{M}$, and the associated ``radiation power'' on node $M$ is henceforth defined as
\begin{equation}
    \label{eq:radiation_power}
    \eta_{M} := T_{M}^4
\end{equation}

The volume integrals in \eqref{eq:SVF_Galerkin} are computed in the usual manner by looping over the elements of $\mathcal{T}_h$ and using standard Gauss quadrature. For the radiation surface integral, we follow the approach presented in the Abaqus manual \cite{Abaqus_2011}, and treat the radiation flux $q_i$ as a constant on each facet $i \in \mathcal{M}_h$, so that this term, which we shall denote by $Q_{N}$, is approximated as
\begin{align}
  \label{eqn:flux_residual}
  Q_{N} :=
  \int_{\Gamma} q(T_h) \varphi_{N}(x) \; \text{d}S  \approx \sum_{i=1}^{\nf} q_i \int_{\Gamma_i} \varphi_{N}(x) \; \text{d}S
\end{align}
In the preceding few equations, we have introduced the notational convention of using $N$, $M$ to index node-based quantities, and $i$, $j$ to index facet-based quantities, and this approach will be continued throughout the rest of this section. To complete the description of the finite element discretization, we need an expression for $q_i$, the radiation flux to the $i$th facet, in terms of all other facets in the radiation cavity. This is discussed in \S\ref{sec:rad_flux_discrete}.


\subsection{Radiation flux discretization\label{sec:rad_flux_discrete}}
The key ingredient in modeling the interaction between discretized facets $i$ and $j$ of an enclosure is the so-called view factor matrix, whose entries are given by \cite{Chapman_1984}
\begin{equation}
  \label{eq:Fdef_new_location}
      F_{ij} = \int_{A_i} \int_{A_j} \frac{\cos{\phi_i}\cos{\phi_j}}{\pi R^2} \, \dx{A_j} \, \dx{A_i},
\end{equation}
\noindent where $A_i$ is the area of facet $i$, $R$ is the distance between facets $i$ and $j$, and $\phi_i$ is the angle between the line segment joining the facet centroids and the plane-normal vector of facet $i$. We note in particular that the diagonal of the view factor matrix consists of all zeros, since \eqref{eq:Fdef_new_location} is only well-defined when $i \neq j$. 

In this work we compute entries of $F$ via numerical quadrature on the facets. The individual entries $F_{ij}$ are non-polynomial and hence may require specialized quadrature rules. A typical approach, which we follow in our implementation, is to apply adaptive Gauss quadrature in which the order of the quadrature rule is increased for facets that are close together (hence strongly interacting), and decreased for facets that are further away from one another. In our testing this leads to accurate and computationally efficient computation of the view factor matrix. However, we note that more sophisticated integration techniques can be applied to the $F_{ij}$ if needed, e.g.\ a detailed discussion for a more efficient treatment of singular integrals representing the interaction of neighboring facets may be found in \cite{Sauter_2011}.

A secondary quantity, which depends on the view factor matrix, is the so-called reflection matrix, which is given by
\begin{equation}\label{eq:reflection_matrix}
    C_{ij} = \delta_{ij} - \frac{1 - \epsilon_i}{A_i} F_{ij} \; \Leftrightarrow \; C = I - \Lambda F,
\end{equation}
where $I$ denotes the identity matrix and $\Lambda$ is a diagonal scaling matrix with $\Lambda_{ii}=(1 - \epsilon_i)/{A_i}$.
With the preceding definitions in mind, the radiation flux to facet $i$ can be written \cite{Abaqus_2011} using the previously defined quantities as
\begin{equation}
   \label{eq:qi_new_location}
   q_i = \frac{\sigma\epsilon_i}{A_i}\sum_j \epsilon_j \sum_k F_{ik} C^{-1}_{kj}(\eta_j - \eta_i)
\end{equation}
where the radiation power $\eta_i$ on facet $i$ is treated as a constant given by the nodal averaging/projection formula
\begin{align}
  \label{eq:facet_power}
  \eta_i := \sum_{M=1}^{\nn} \frac{1}{A_i} \int_{\Gamma_i} \eta_{M} \varphi_{M}(x) \; \text{d}S
\end{align}
where $A_i$ is the area of facet $i$, whose domain is $\Gamma_i$. We observe that the sum in \eqref{eq:facet_power} nominally goes over all nodes $M$ in the finite element mesh, but in fact only a small number of nodal basis functions $\varphi_{M}$ are non-zero on any given facet $\Gamma_i$ due to the localized nature of the finite element basis functions. Therefore, many of the terms in the sum are zero. Associated with the formula in \eqref{eq:facet_power}, it is illustrative to define the ``projection'' operator
\begin{align}
  \label{eqn:projection_operator}
  P_{iM} := \frac{1}{A_i} \int_{\Gamma_i} \varphi_{M}(x) \; \text{d}S
\end{align}
so that
\begin{align}
  \label{eqn:P_fi_definition}
  \eta_i = \sum_{M=1}^{\nn} P_{iM} \eta_{M}
\end{align}
The projection operator $P_{iM}$ therefore maps node-based quantities to facet-based quantities, and we can think of it as a sparse, rectangular $\nf \times \nn$ matrix for simplicity.




Combining all of the preceding ingredients, we finally get the following expression for $Q_{N}$, originally defined in \eqref{eqn:flux_residual} (note: all sums go from 1 to $\nf$ unless otherwise noted)
\begin{align}
  \nonumber
  Q_{N} &:= \int_{\Gamma} q(T_h) \varphi_{N}(x) \; \text{d}S
  \\
  \nonumber
  &\approx \sum_{i} q_i \int_{\Gamma_i} \varphi_{N}(x) \; \text{d}S
  \\
  \nonumber
  &= \sum_{i} \frac{\sigma\epsilon_i}{A_i}\sum_j \epsilon_j \sum_k F_{ik} C^{-1}_{kj}(\eta_j - \eta_i) \int_{\Gamma_i} \varphi_{N}(x) \; \text{d}S
  \\
  \nonumber
  &= \sum_{i} \left(\sum_j R_{ij} \eta_j - \sum_k R_{ik} \eta_i\right) \left( \frac{1}{A_i}\int_{\Gamma_i} \varphi_{N}(x) \; \text{d}S \right)
  \\
  \nonumber
  &= \sum_{i} \sum_j \bar{R}_{ij} \eta_j P_{iN}
  \\
  \nonumber
  &= \sum_{i} \sum_j \bar{R}_{ij} 
  \left(\sum_{M=1}^{\nn} P_{jM} \eta_{M} \right)   
  P_{iN}
  \\
  \label{eqn:flux_residual_2}
  &= \sum_{i} \sum_j \sum_{M=1}^{\nn} P_{iN} \bar{R}_{ij} P_{jM} \eta_{M}  
\end{align}
where we have used \eqref{eqn:projection_operator}, \eqref{eqn:P_fi_definition} and, for simplicity of notation, we have defined the tensors $R$ and $\bar{R}$ as
\begin{align}
\label{eq:Rdef}
R_{ij} &:= \sigma \epsilon_i \epsilon_j \sum_k F_{ik} C^{-1}_{kj} 
\\
\label{eq:Rbardef}
\bar{R}_{ij} &:= R_{ij} - \delta_{ij}\sum_{k} R_{ik}
\end{align}
We can write the double facet sum from \eqref{eqn:flux_residual_2} as
\begin{align}
    \label{eqn:S_def}
    S_{N M} := \sum_{i} \sum_j P_{iN} \bar{R}_{ij} P_{jM} 
\end{align}
or, employing matrix notation
\begin{align}
\label{eqn:S_def_matrix}
    S := P^\top \bar{R} P
\end{align}
which gives us the following final expression for $Q_{N}$
\begin{align}
\label{eqn:flux_residual_3}
    Q_{N} = \sum_{M=1}^{\nn} S_{N M} \eta_{M}
\end{align}

\subsection{Time discretization and Newton iteration\label{sec:time_discretization}}
In the time domain, we discretize \eqref{eq:SVF_Galerkin} using a standard implicit Euler scheme. Let $T^{n}_h$ be the approximate solution at time level $t_n$, and consider a fixed time step of size $\Delta t$ such that $t_{n+1} = t_n + \Delta t$. The $N$th component of the discrete residual vector associated with \eqref{eq:SVF_Galerkin} is then given by
\begin{align}
    \nonumber
    \mathscr{R}_{N}(T^{n+1}_h) &:= \int_{\Omega} \rho c_p \left(\frac{ T^{n+1}_h - T^n_h}{\Delta t} \right) \varphi_{N} \, \dx{V} 
      + \int_{\Omega} k \nabla T^{n+1}_h \cdot \nabla \varphi_{N} \, \dx{V} 
      \\
      \label{eq:SVF_Galerkin_time}
      &- Q^{n+1}_N - \int_{\Omega} f \varphi_{N} \, \dx{V} 
\end{align}
The solution of the nonlinear system of equations $\mathscr{R} = 0$ implied by \eqref{eq:SVF_Galerkin_time} is performed at each time step via Newton--Krylov \cite{Dembo_1982,Knoll_2004} iterations, which require the assembly of the Jacobian $J$ associated with $\mathscr{R}$ based on the current temperature iterate.

Note that the Jacobian contribution due to the volume integral terms, which we shall denote by $\Jsparse$, will be sparse due to the usual properties of the Galerkin finite element method. The flux term $Q_N^{n+1}$, on the other hand, is both nonlinear and non-local, as discussed previously. The Jacobian contribution for this term, $\Jcav$, is thus a dense matrix given by
\begin{equation}\label{eq:Jcav}
    \Jcav_{NM} := \frac{\partial Q_N}{\partial T_M} = \sum_{K=1}^{\nn} {S}_{NK} \frac{\partial }{\partial T_M} \left(T^4_K\right) = 4{S}_{NM} T^3_M
\end{equation}
\noindent where we have dropped the $n+1$ superscript for clarity. We also emphasize that the matrix $S$, as defined in (\ref{eqn:S_def_matrix}), forms the core part of $\Jcav$, and also is temperature independent. The total Jacobian
\begin{equation}\label{eq:Jtotal}
    J = \Jsparse + \Jcav
\end{equation}
therefore has both sparse and dense contributions, and requires specialized handling, which is described in Section \ref{sec:iteration}.

\subsection{Iterative approach\label{sec:iteration}}

In case of a large number of facets, the matrix size of the (dense) cavity part $\Jcav$ of the Jacobian becomes prohibitive, making direct solution methods for the solution of the linear system in each nonlinear iteration impractical. We hence resort to preconditioned iterative methods (such as GMRES), requiring the implementation of the action of the Jacobian and a suitable preconditioner onto an arbitrary vector $x$.

Thanks to the additive decomposition of the total Jacobian \eqref{eq:Jtotal}, we may consider the two terms in the product
\begin{equation}\label{eq:Jproduct}
    Jx = \Jsparse x + \Jcav x
\end{equation}
in isolation. Since it is inexpensive to compute, we assume that $\Jsparse$ is available in fully assembled (sparse) form, and the resulting product $\Jsparse x$ can therefore be computed efficiently by standard matrix-vector operations. In addition to being conveniently fully assembled, we have also observed that $\Jsparse$ is suitable for use as a preconditioner. Our preconditioned iteration hence reads
\begin{equation}
    \label{eq:Jprec}
    x \mapsto (\Jsparse)^{-1} \left(\Jsparse x + \Jcav x\right)
\end{equation}
where $(\Jsparse)^{-1}$ is computed explicitly via a sparse LU decomposition. We note that more refined choices of preconditioners are possible (in particular involving additional terms for the cavity part of the Jacobian), but we have not investigated these in more detail since, in our examples, we generally observed good performance (small numbers of GMRES iterations) with this simple choice.

The approach for efficiently storing and computing matrix-vector products with $\Jcav$ relies on a block low-rank approximation of the view factor matrix $F$, and a subsequent approximation of the inverse of the reflection matrix $C$; these are described in detail in Sections \ref{sec:view_factor} and \ref{sec:cavity_jac_assembly}. If the material properties (thermal conductivity, specific heat, etc.)\ are temperature- and time-independent, then $\Jsparse$ is itself also independent of the temperature, and thus theoretically only needs to be computed once and reused for all subsequent iterations. Our present implementation is designed to handle the general case of temperature-dependent material properties, and thus does not currently take advantage of this particular optimization. This choice is justified by the solver performance results reported in Section \ref{sec:numerical_experiments}, which show that, especially in the fine grid limit, the time spent recomputing $\Jsparse$ at each iteration is small compared to the time spent in the computation of $\Jcav$.

\section{Block low-rank approximation of the view factor matrix\label{sec:view_factor}}
The view factor matrix $F$ defined in \eqref{eq:Fdef_new_location} represents the radiation interaction between different facets from the mesh $\mathcal{M}_h$.  The particular form of \eqref{eq:Fdef_new_location} immediately gives rise to the following two observations (cf. \cite{Potter2022}):
\begin{enumerate}[(i)]
    \item Although we do not consider occlusion detection in the current work, matrix blocks corresponding to occluded groups of facets are identically zero. This may result in major computational gains, as the corresponding matrix entries need not be computed or stored explicitly.
    \item\label{item:R2} The radiation interaction decays with increasing distance $R$. Not only does this make the entries $F_{ij}$ smaller in absolute value, but it also leads to a weak interaction between groups of facets that are far apart. This enables data-sparse approximations of the corresponding matrix blocks of the view factor matrix.
\end{enumerate}

Our focus in this paper is to explore aspect \ref{item:R2} using hierarchical block low-rank approximations of the view factor matrix. The key elements of this approximation are:
\begin{enumerate}[(a)]
    \item\label{item:geo_clustering} Geometric clustering of the facets via a $k$-d tree to hierarchically subdivide the row and column index set $\mathcal{I}=\{1,...,n\}$ of $F$.
    \item\label{item:subdivision} Using the subdivision from \ref{item:geo_clustering}, the block index set $\mathcal{I}\times\mathcal{I}$ is hierarchically subdivided into smaller blocks until the corresponding facet groups are sufficiently far apart. To prevent unreasonably small matrix blocks, the subdivision is also stopped if a user-defined minimal block size has been reached.
    \item\label{item:distant} Matrix blocks corresponding to distant groups of facets are approximated by low-rank matrices. The rest of the blocks are stored as dense matrices.
\end{enumerate}
In the following, we provide a more detailed description of the techniques used for \ref{item:geo_clustering}, \ref{item:subdivision} and \ref{item:distant}.

\subsection{Geometric Clustering\label{sec:geometry}}
To group facets from the mesh $\mathcal{M}_h$ into clusters, we use a $k$-d tree based on the centroids of the individual facets. We note that an alternative, yet similar approach would be to consider a clustering approach based on the bounding boxes of the facets.

With the index set $\mathcal{I} := \{1,...,n\}$ corresponding to the centroids of all facets as a starting point, the clustering proceeds iteratively by geometrically splitting the point cloud of centroids via a hyperplane. Once the split has been defined, the iteration continues recursively for the two disjoint sets of centroids. The procedure for an arbitrary index set $\tau\subset\mathcal{I}$ is defined in Algorithm \ref{alg:kdtree}. 

\begin{algorithm}[!ht]
\caption{\texttt{BuildIndexTree}($\tau \subset \mathcal{I}$)\label{alg:kdtree}}
\begin{algorithmic}
\If{ $\#\tau \leq n_\text{min}$}
    \State $T \gets \{\tau\}$
    \Comment{leaf index set}
\Else 
    \State $\tau =: \tau_\text{left} \mathbin{\mathaccent\cdot\cup} \tau_\text{right}$
    \Comment{split using a hyperplane}
    \State $T_\text{left} \gets$ \texttt{BuildIndexTree}$(\tau_\text{left})$
    \State $T_\text{right} \gets$ \texttt{BuildIndexTree}$(\tau_\text{right})$
    \State $T \gets T_\text{left} \uplus T_\text{right}$
    \Comment{left and right subtrees}
\EndIf
\end{algorithmic}
\end{algorithm}

The initial call to Algorithm \ref{alg:kdtree} reads
\begin{equation}
    T_\mathcal{I} := \texttt{BuildIndexTree}(\mathcal{I}),
\end{equation}
where the obtained \emph{index tree} $T_\mathcal{I}$ represents the hierarchical split of all centroids. We note that in Algorithm \ref{alg:kdtree} the recursion is only continued in case the cardinality of the index set $\tau$ is larger than a user-defined minimal leaf size $n_\text{min}\geq 1$ (say $n_\text{min} \sim 100$). This is to prevent the creation of arbitrarily small matrix blocks which would negatively impact the performance of the overall hierarchical approach.

Since the splitting in Algorithm \ref{alg:kdtree} is unaware of the original ordering of the facets in $\mathcal{M}_h$, the nodes $\tau\in T_\mathcal{I}$ generally correspond to non-contiguous index sets. To associate the obtained splitting to contiguous row and column index sets of the view factor matrix $F$, it is hence convenient to permute the rows and columns of $F$ accordingly.

\subsection{Block Subdivision}
Using the hierarchical subdivision of the index set $\mathcal{I}$, we next introduce a block subdivision of $\mathcal{I}\times \mathcal{I}$ defining the eventual matrix blocks of the view factor matrix $F$. As before, we proceed by recursively subdividing the block index set $\mathcal{I}\times \mathcal{I}$ into smaller subblocks until a stopping criterion is fulfilled.

Since we expect that groups of facets which are far apart interact only weakly, we use a standard geometrically-based \emph{admissibility condition} (cf. \cite{Hackbusch2015}) as a stopping criterion. Given two index sets $\sigma, \tau \subset \mathcal{I}$, we denote the corresponding matrix block $\sigma \times \tau$ as \emph{admissible} if
\begin{equation}
    \min(\operatorname{diam}(A_\sigma),\operatorname{diam}(A_\tau)) \leq c \operatorname{dist}(A_\sigma,A_\tau),
\end{equation}
where $A_\sigma, A_\tau$ denote the point clouds of centroids of the facets belonging to index sets $\sigma, \tau$, respectively. The recursive procedure for the block subdivision is defined in Algorithm \ref{alg:blocktree}. 

\begin{algorithm}[!ht]
\caption{\texttt{BuildBlockTree}($\sigma,\tau\in T_\mathcal{I}$)}\label{alg:blocktree}

\begin{algorithmic}
\If{the block $\sigma \times \tau$ is admissible}
    \State $T \gets \{\sigma\times\tau\}$
    \Comment{leaf block}
\Else 
    \State $T \gets \emptyset$
    \For{$\sigma^\prime \in \operatorname{child}(\sigma)$ and $\tau^\prime \in \operatorname{child}(\tau)$}
        \Comment{recursion to children}
        \State $T_\text{sub} \gets\texttt{BuildBlockTree}(\sigma^\prime,\tau^\prime$)
        \State $T \gets T\uplus T_\text{sub}$
    \EndFor
\EndIf
\end{algorithmic}
\end{algorithm}

The initial call to Algorithm \ref{alg:blocktree} reads
\begin{equation}
    T_{\mathcal{I}\times\mathcal{I}} := \texttt{BuildBlockTree}(\mathcal{I},\mathcal{I}),
\end{equation}
where the obtained \emph{block tree} $T_{\mathcal{I}\times\mathcal{I}}$ represents the hierarchical block subdivision of the view factor matrix $F$.

\begin{example}
Consider a simple equidistant arrangement of the centroids of $n$ facets along a straight line. For a trivial minimal leaf size of $n_\text{min}=1$, Algorithms \ref{alg:kdtree} and \ref{alg:blocktree} would result in a block structure as depicted in Figure \ref{fig:blocktree}. 

\input{images/hierarchicalmatrix}

\end{example}

\begin{example}

Consider two finite parallel flat plates of size $L \times L$, separated by a distance $L$, and let their surfaces be successively refined so that in total there are $40\times 40$, $64\times 64$, and $106\times 106$ facets, respectively, on each surface as shown in Fig.~\ref{fig:flat_plates}. Then, using a leaf size of 128, we obtain the block low-rank structure shown in Fig.~\ref{fig:Franklin_test_case_block_structure_all}. Note that the ``on-diagonal'' blocks are initially represented as dense matrices on the coarse grids, but they are eventually revealed to have some internal low-rank structure as the mesh is refined.

\begin{figure}[b]
\centering
\includegraphics[width=0.5\textwidth]{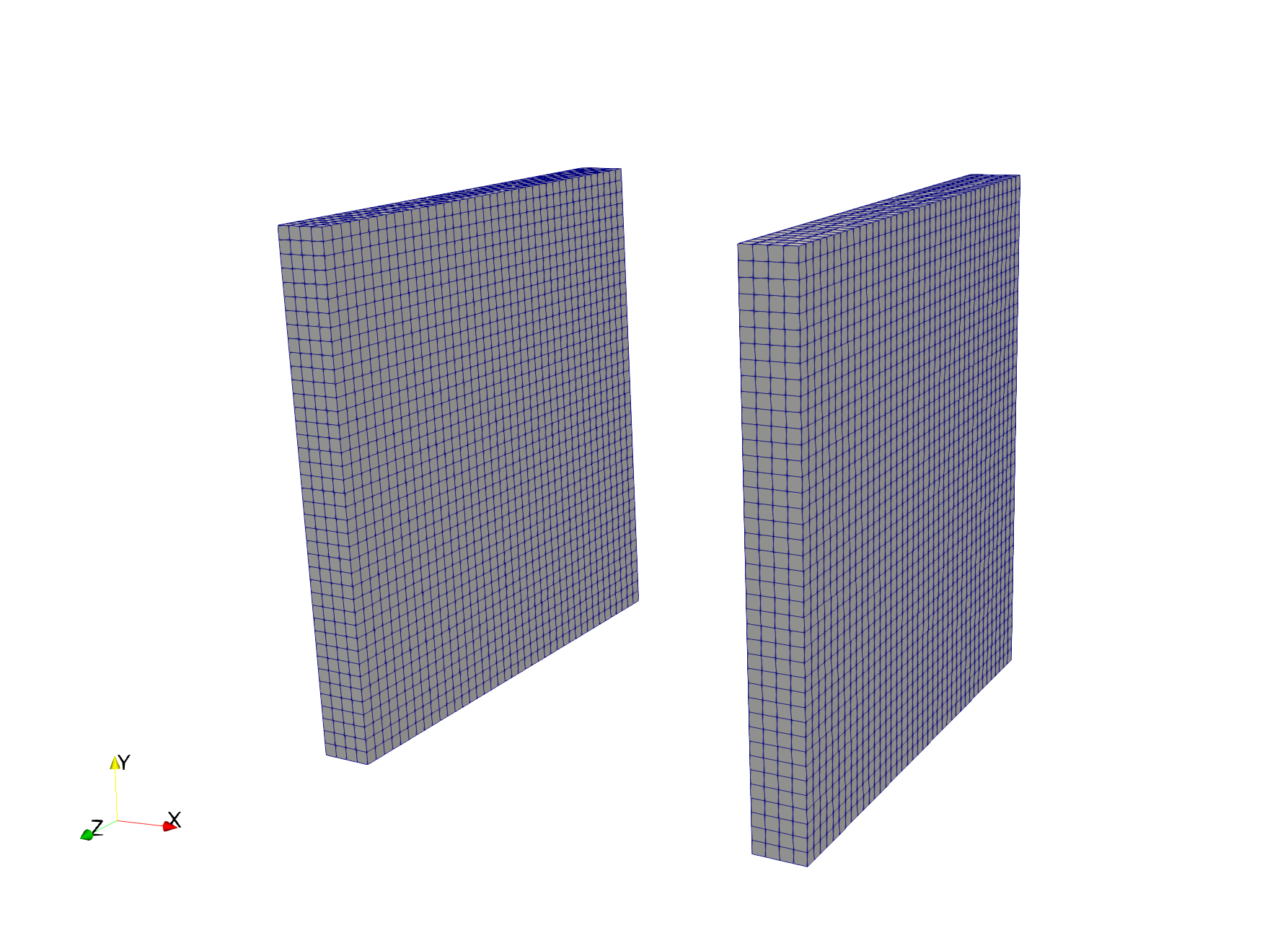}
\caption{Separated parallel flat plates with 3200 facets on the cavity surface.\label{fig:flat_plates}}
\end{figure}

\begin{figure}[!htb]
\centering
\begin{subfigure}[t]{.5\textwidth}
  \centering
  \includegraphics[width=\linewidth]{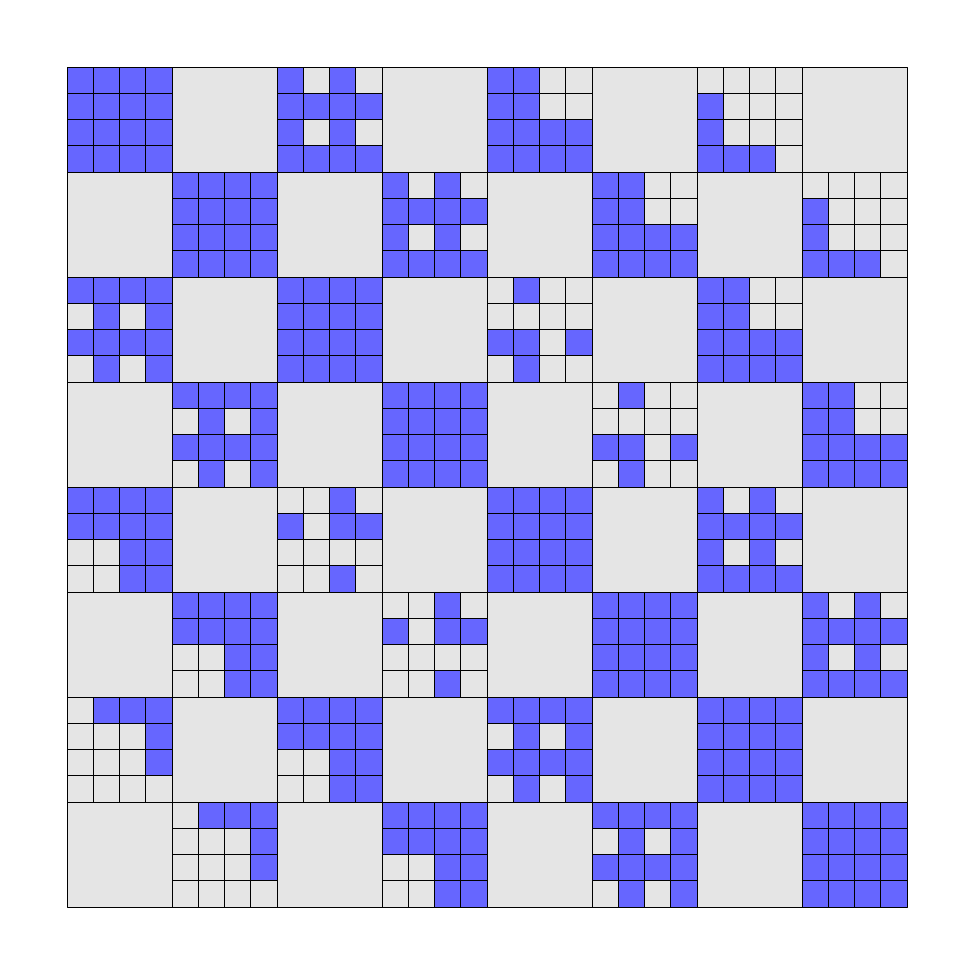}
  \caption{\num{3200} facets \label{fig:Franklin_test_case_block_structure}}
 \end{subfigure}%
\begin{subfigure}[t]{.5\textwidth}
  \centering
  \includegraphics[width=\linewidth]{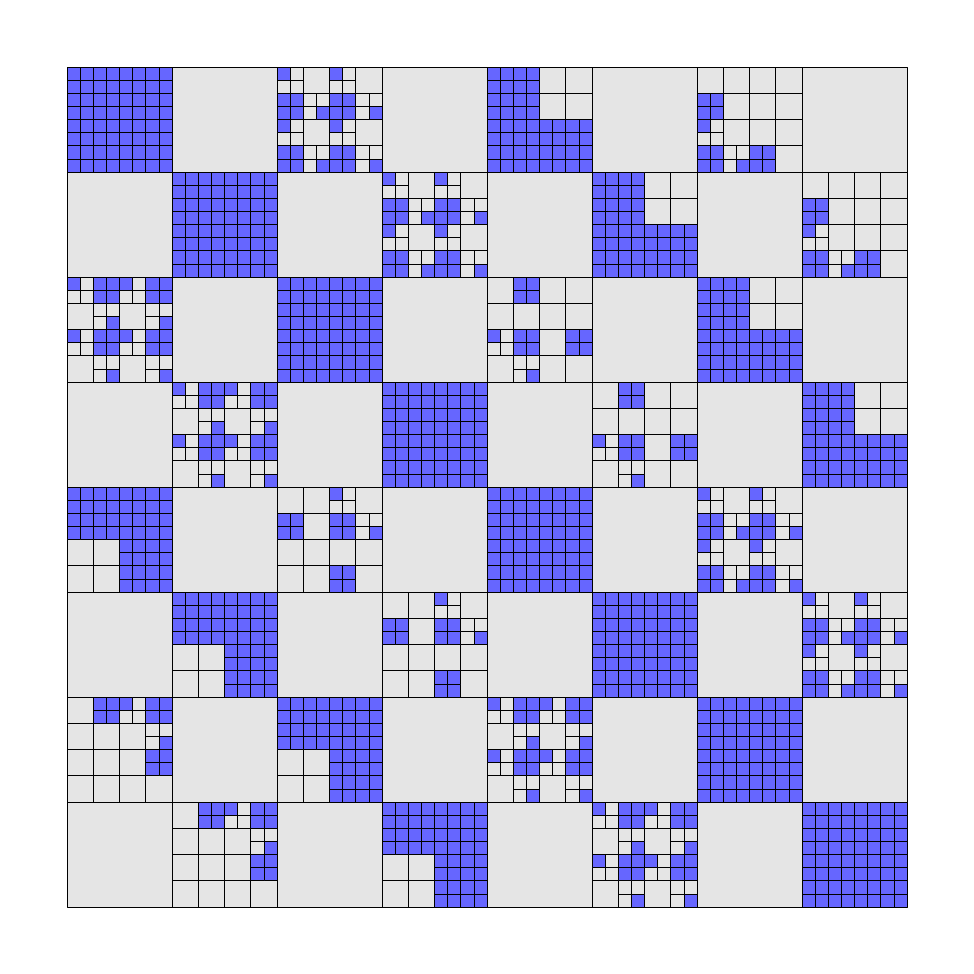}
  \caption{\num{8192} facets \label{fig:Franklin_test_case_8k_block_structure}}
\end{subfigure}
\\
\begin{subfigure}[t]{.5\textwidth}
  \centering
  \includegraphics[width=\linewidth]{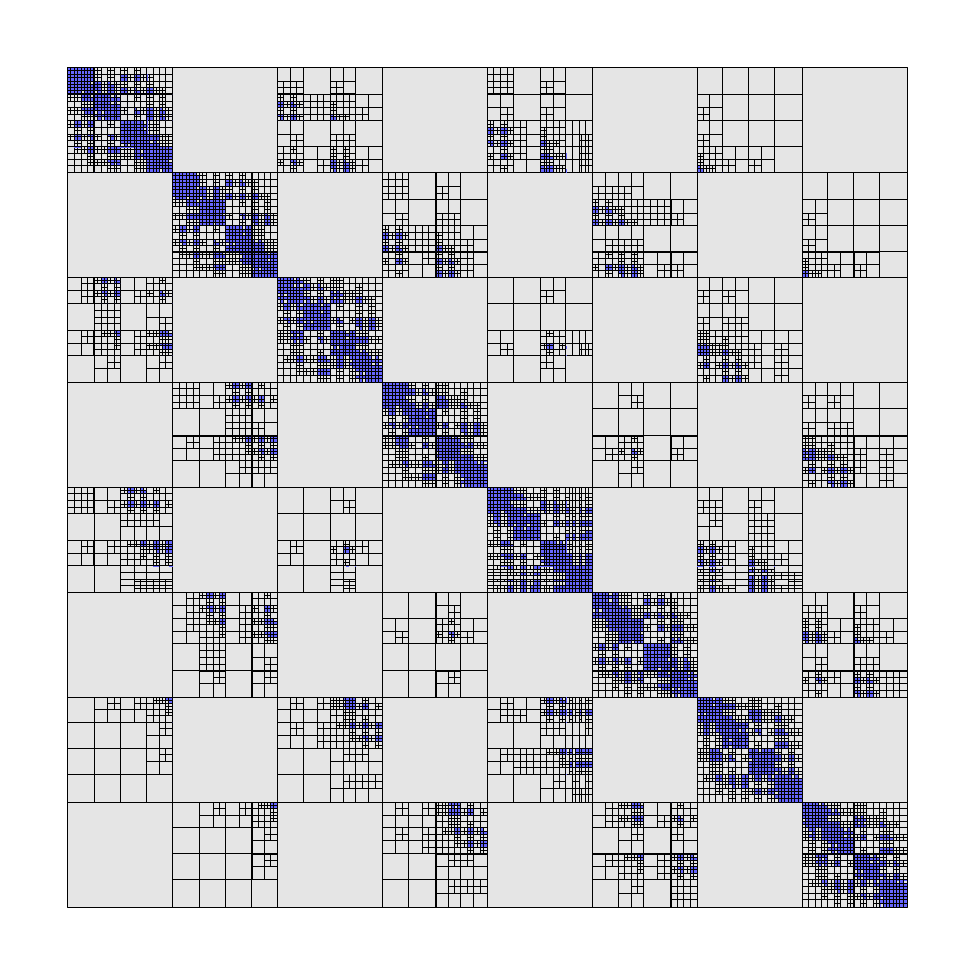}
  \caption{\num{22472} facets \label{fig:Franklin_test_case_22k_block_structure}}
\end{subfigure}

\caption{View factor matrix block low-rank approximation for successively refined parallel flat plates. These block structures were computed numerically using the approach described in this Section.\label{fig:Franklin_test_case_block_structure_all}}
\end{figure}

\end{example}

\subsection{Adaptive Cross Approximation}




For the low-rank approximation of data sparse blocks one can use the singular value decomposition (SVD) to compute the low-rank decomposition, but it is well know that this becomes computationally expensive for large blocks. An alternative approach that is typically more computationally efficient is to apply adaptive cross approximation (ACA, cf. \cite{Bebendorf2000}). Using ACA we can construct a $k$-rank approximation of a data sparse matrix $A \in \mathbb{R}^{m, n}$ just in $\mathcal{O}(m n k)$ time, which is significantly faster than classical SVD.

\begin{figure}[t]
\centering
\begin{tikzpicture}
    \tikzstyle{Sfill} = [fill=blue!60];
    \tikzstyle{Rfill} = [fill=black!30];
    \coordinate (X) at (0.15,1.5);  
    \coordinate (Y) at (1.5,0.15);
    \coordinate (XY) at (0.15,0.15);
    \coordinate (P1) at ($(-1.0,0.0)-10*(1,0 -|X)$);
    \coordinate (Q1) at ($10*(1,0 -|X) + 10*(0,1 |- Y)$);
    \coordinate (C1) at ($(P1)-10*(1,0 -|X)$);
    \coordinate (C2) at ($(P1)-6*(1,0 -|X)$);  
    \coordinate (C3) at ($(P1)-0*(1,0 -|X)$);    
    \coordinate (R1) at ($(P1)+7*(0,1 |-Y)$);
    \coordinate (R2) at ($(P1)+3*(0,1 |-Y)$);
    \coordinate (R3) at ($(P1)-3*(0,1 |-Y)$);    
    \coordinate (A11) at ($(1,0 -| C1) + (0,1 |- R1)$);
    \coordinate (A12) at ($(1,0 -| C1) + (0,1 |- R2)$);  
    \coordinate (A13) at ($(1,0 -| C1) + (0,1 |- R3)$);  
    \coordinate (A21) at ($(1,0 -| C2) + (0,1 |- R1)$);  
    \coordinate (A22) at ($(1,0 -| C2) + (0,1 |- R2)$);  
    \coordinate (A23) at ($(1,0 -| C2) + (0,1 |- R3)$);  
    \coordinate (A31) at ($(1,0 -| C3) + (0,1 |- R1)$);  
    \coordinate (A32) at ($(1,0 -| C3) + (0,1 |- R2)$);  
    \coordinate (A33) at ($(1,0 -| C3) + (0,1 |- R3)$);   
    \coordinate (D1) at (1.1,0.0);
    \coordinate (D2) at ($(D1)+2*(1,0 -|X)$);
    \coordinate (D3) at ($(D1)+4*(1,0 -|X)$);
    \coordinate (S1) at (5.0,1.35);  
    \coordinate (S2) at ($(S1)-2*(0,1 |- Y)$);
    \coordinate (S3) at ($(S1)-4*(0,1 |- Y)$);
    \coordinate (B1) at ($0.5*(1,0 -| D3) + 0.5*(1,0 -| X) + 0.5*(1,0 -| S1) - 0.5*(1,0 -| Y) 
                         +(0,1 |- X) - 3*(0,1 |- Y)$);
    \coordinate (B11) at ($(B1) - 2*(1,0 -| X) + 2*(0,1 |- Y)$);
    \coordinate (B12) at ($(B1) + 0*(1,0 -| X) + 2*(0,1 |- Y)$);
    \coordinate (B13) at ($(B1) + 2*(1,0 -| X) + 2*(0,1 |- Y)$);
    \coordinate (B21) at ($(B1) - 2*(1,0 -| X) + 0*(0,1 |- Y)$);
    \coordinate (B22) at ($(B1) - 0*(1,0 -| X) + 0*(0,1 |- Y)$);
    \coordinate (B23) at ($(B1) + 2*(1,0 -| X) + 0*(0,1 |- Y)$);    
    \coordinate (B31) at ($(B1) - 2*(1,0 -| X) - 2*(0,1 |- Y)$);
    \coordinate (B32) at ($(B1) + 0*(1,0 -| X) - 2*(0,1 |- Y)$);
    \coordinate (B33) at ($(B1) + 2*(1,0 -| X) - 2*(0,1 |- Y)$);
    \filldraw[fill=black!10] ($(P1)-(Q1)$) rectangle ($(P1)+(Q1)$);
    \filldraw[Rfill] ($(C1)-(X)$) rectangle ($(C1)+(X)$);
    \filldraw[Rfill] ($(C2)-(X)$) rectangle ($(C2)+(X)$);
    \filldraw[Rfill] ($(C3)-(X)$) rectangle ($(C3)+(X)$);
    \filldraw[Rfill] ($(R1)-(Y)$) rectangle ($(R1)+(Y)$);
    \filldraw[Rfill] ($(R2)-(Y)$) rectangle ($(R2)+(Y)$);
    \filldraw[Rfill] ($(R3)-(Y)$) rectangle ($(R3)+(Y)$);    
    \filldraw[Sfill] ($(A11)-(XY)$) rectangle ($(A11)+(XY)$);
    \filldraw[Sfill] ($(A12)-(XY)$) rectangle ($(A12)+(XY)$);
    \filldraw[Sfill] ($(A13)-(XY)$) rectangle ($(A13)+(XY)$);
    \filldraw[Sfill] ($(A21)-(XY)$) rectangle ($(A21)+(XY)$);
    \filldraw[Sfill] ($(A22)-(XY)$) rectangle ($(A22)+(XY)$);
    \filldraw[Sfill] ($(A23)-(XY)$) rectangle ($(A23)+(XY)$);
    \filldraw[Sfill] ($(A31)-(XY)$) rectangle ($(A31)+(XY)$);
    \filldraw[Sfill] ($(A32)-(XY)$) rectangle ($(A32)+(XY)$);
    \filldraw[Sfill] ($(A33)-(XY)$) rectangle ($(A33)+(XY)$);
    \filldraw[Rfill] ($(D1)-(X)$) rectangle ($(D1)+(X)$);
    \filldraw[Rfill] ($(D2)-(X)$) rectangle ($(D2)+(X)$);    
    \filldraw[Rfill] ($(D3)-(X)$) rectangle ($(D3)+(X)$); 
    \filldraw[Sfill] ($(B11)-(XY)$) rectangle ($(B11)+(XY)$); 
    \filldraw[Sfill] ($(B12)-(XY)$) rectangle ($(B12)+(XY)$); 
    \filldraw[Sfill] ($(B13)-(XY)$) rectangle ($(B13)+(XY)$); 
    \filldraw[Sfill] ($(B21)-(XY)$) rectangle ($(B21)+(XY)$); 
    \filldraw[Sfill] ($(B22)-(XY)$) rectangle ($(B22)+(XY)$); 
    \filldraw[Sfill] ($(B23)-(XY)$) rectangle ($(B23)+(XY)$); 
    \filldraw[Sfill] ($(B31)-(XY)$) rectangle ($(B31)+(XY)$); 
    \filldraw[Sfill] ($(B32)-(XY)$) rectangle ($(B32)+(XY)$); 
    \filldraw[Sfill] ($(B33)-(XY)$) rectangle ($(B33)+(XY)$); 
    \filldraw[Rfill] ($(S1)-(Y)$) rectangle ($(S1)+(Y)$);
    \filldraw[Rfill] ($(S2)-(Y)$) rectangle ($(S2)+(Y)$);
    \filldraw[Rfill] ($(S3)-(Y)$) rectangle ($(S3)+(Y)$);          
    \draw ($(C1)+(0,1 |- X) $) node[anchor=south] {1};
    \draw ($(C2)+(0,1 |- X) $) node[anchor=south] {3};
    \draw ($(C3)+(0,1 |- X) $) node[anchor=south] {6};
    \draw ($(-0.15,0.0)+(R1)-(1,0 -| Y) $) node[anchor=east]  {2};
    \draw ($(-0.15,0.0)+(R2)-(1,0 -| Y) $) node[anchor=east]  {4};
    \draw ($(-0.15,0.0)+(R3)-(1,0 -| Y) $) node[anchor=east]  {7};
    \draw ($(D1)+(0,1 |- X) $) node[anchor=south] {1};
    \draw ($(D2)+(0,1 |- X) $) node[anchor=south] {3};
    \draw ($(D3)+(0,1 |- X) $) node[anchor=south] {6};
    \draw ($(S1)+(1,0 -| Y) $) node[anchor=west]  {2};
    \draw ($(S2)+(1,0 -| Y) $) node[anchor=west]  {4};
    \draw ($(S3)+(1,0 -| Y) $) node[anchor=west]  {7};
    \draw ($(0.0,0.0)$) node {$\approx$};
\end{tikzpicture}
\caption{Low rank approximation by Adaptive Cross Approximation (ACA) with row and column indices $\mathcal{R}=\{2,4,7\}$ and $\mathcal{C}=\{1,3,6\}$, respectively}
\label{fig:ACA}
\end{figure}
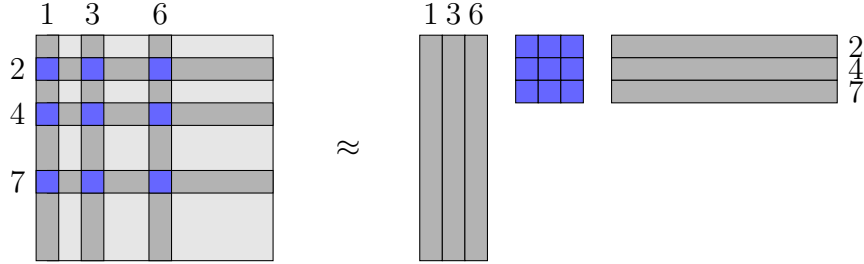

The main idea of ACA is to iteratively determine sets $\mathcal{R} = \{i_1,\ldots,i_k\}$, $\mathcal{C} = \{j_1,\ldots,j_k\}$ of row and column indices, respectively, to construct a low-rank approximation of admissible matrix blocks $X$ of $F$ of the form
\begin{equation}\label{eq:aca}
    X \approx X(:,\mathcal{C}) \left(X(\mathcal{R},\mathcal{C})\right)^{-1} X(\mathcal{R}, :),
\end{equation}
see Figure \ref{fig:ACA} for an illustration. The representation \eqref{eq:aca} may then be further postprocessed (or truncated) to obtain low-rank factors $X \approx \tilde{X} := UV^{\top}$ in standard form.

There are multiple approaches to how the pivoting should be performed (some of them provide accuracy guarantees related to best low-rank approximation), but exact details go beyond the scope of this paper. A more detailed discussion on the different techniques can be found in \cite{Borm_2006}. In Section \ref{sec:numerical_experiments} we investigate different approximation tolerances $\epsrel$ which can be set for ACA to satisfy:

\begin{equation}\label{eq:erel}
    \frac{\| X - \tilde{X} \|}{\| X \|} \leq \epsrel
\end{equation}

\noindent for all blocks $X$ as defined above. We will analyze the resulting accuracy and demonstrate the computational advantage of ACA for practical cavity radiation problems with respect to the choice of $\epsrel$ in the solver.

\section{Cavity Jacobian assembly\label{sec:cavity_jac_assembly}}
The iterative approach introduced in Section \ref{sec:iteration} requires us to efficiently construct the cavity part $\Jcav$ of the Jacobian from \eqref{eq:Jcav} and to implement its action on an arbitrary vector. We recall from \eqref{eqn:S_def_matrix} that the matrix-vector product with $\Jcav$ is defined through a (sparse) projection operator $P$ and a dense coupling matrix $\bar{R}$ which depends on the view factor matrix $F$ and the inverse of the reflection matrix $C$ via \eqref{eq:Rdef}, \eqref{eq:Rbardef}.

We now assume that the view factor matrix $F$ has been constructed by the techniques introduced in Section \ref{sec:view_factor} and is hence available in a hierarchical block low-rank format. It is then straightforward to obtain the reflection matrix $C = I-\Lambda F$ from \eqref{eq:reflection_matrix} within the same hierarchical format as well, since the diagonal matrix $\Lambda$ only introduces a simple row scaling in each matrix block of $F$, and the identity matrix $I$ can be added to the diagonal (dense) blocks without modification of the block structure.

The key computational task is hence to efficiently compute an (approximate) factorization of $C$. For simplicity, we construct this factorization within the same hierarchical block low-rank format as $F$, though alternative choices are possible (cf. \cite{Hackbusch2015}). The remaining task within our iterative framework is then to implement matrix-vector multiplication with matrices in hierarchical block low-rank format and arbitrary vectors.

\subsection{Hierarchical block LU decomposition of the reflection matrix}
The factorization of the reflection matrix $C \approx LU$ proceeds recursively using the block structure defined through Algorithm \ref{alg:blocktree}. The first block subdivision leads to a block system
\begin{equation*}
    \begin{bmatrix}
        C_{11} & C_{12} \\
        C_{21} & C_{22} 
    \end{bmatrix}
    = 
    \begin{bmatrix}
        L_{11} & 0 \\
        L_{21} & L_{22}
    \end{bmatrix}
    \cdot
    \begin{bmatrix}
        U_{11} & U_{12} \\
        0 & U_{22}
    \end{bmatrix}
\end{equation*}
which can be solved via the following four steps:
\begin{enumerate}
\item\label{item:LU11} perform the LU decomposition $L_{11}U_{11} = C_{11}$,
\item\label{item:frwdsubs} solve for $U_{12} = L_{11}^{-1}C_{12}$ by forward substitution,
\item\label{item:backsubs} solve for $L_{21} = C_{21}U_{11}^{-1}$ by backward substitution,
\item\label{item:LU22} perform the LU decomposition $L_{22}U_{22} = C_{22} - L_{21}U_{12}$.
\end{enumerate}

The LU decompositions in Steps \ref{item:LU11} and \ref{item:LU22} are computed recursively in case the corresponding blocks are further subdivided, whereas dense LU factorizations are employed in leaf blocks. The forward and backward substitutions in Steps \ref{item:frwdsubs} and \ref{item:backsubs} are again implemented by recursive algorithms involving block matrix-matrix products and additions, cf. \cite{Hackbusch2015} for a detailed outline.

It is important to note that matrix additions and multiplications will formally lead to a rank increase in the corresponding low-rank blocks and hence need to be coupled with a suitable rank truncation policy. In our experiments, we have chosen to fix a relative accuracy at the value $\epsrel$ (equal to the approximation accuracy used for the view factor matrix $F$) and to determine the ranks adaptively based on the singular value decay within the low-rank blocks. As a consequence of the truncation, the computed LU decomposition $C \approx LU$ is only approximate, but within a prescribed, controllable accuracy.

\subsection{Matrix-vector assembly}
Once the initial factorization of the reflection matrix has been obtained, it can be re-used within all iterations of the cavity Jacobian assembly. In particular, the action $x \mapsto C^{-1}x$ is approximated by $x \mapsto U^{-1} (L^{-1}x)$ which is implemented as a subsequent forward and backward substitution. Similarly, the assembly of the action $x \mapsto Fx$ is performed recursively, following the block structure of $F$. Most remarkably, the matrix-vector assembly does not involve any further approximation but is performed ``exactly" (within the numerical accuracy of floating point operations).

\section{Numerical experiments\label{sec:numerical_experiments}}

In this section, we consider radiative heat transfer between a collection of spheres arranged in a Fibonacci spiral pattern \cite{Dunlap_1997}. The spiral configuration is chosen because it leads to a distribution of both nearby and distant facets that provide a good illustration of the capabilities of the block low-rank approach presented in the preceding sections. (We also note that there is little occlusion between the various spheres when they are arranged in this manner, which justifies our omission of occlusion detection for this test case.) Since we include the entire exterior boundary of each sphere in the radiation cavity, there will be some facets (which we term ``isolated'' facets) that cannot ``see'' any other facets within the cavity, and thus produce identically zero rows in the view factor matrix.

\subsection{Problem setup and representative results\label{sec:problem_setup}}
\begin{figure}[htb]
\centering
\begin{subfigure}[t]{.48\textwidth}
  \centering
  \includegraphics[width=.95\linewidth]{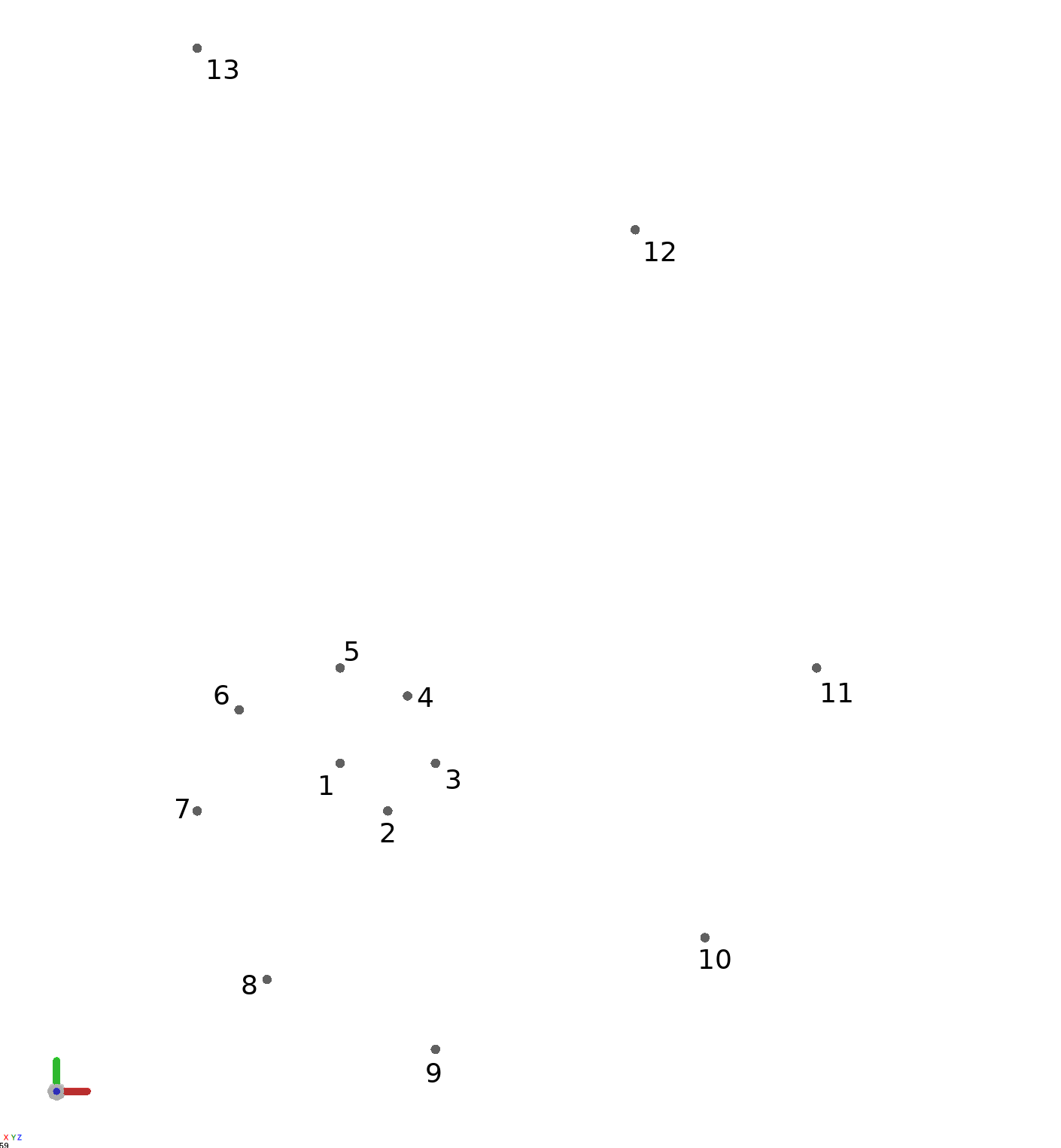}
  \caption{Spiral configuration \label{fig:spiral_overhead_view_labels}}
 \end{subfigure}%
\begin{subfigure}[t]{.48\textwidth}
  \centering
  \raisebox{.25\linewidth}{\includegraphics[width=.95\linewidth]{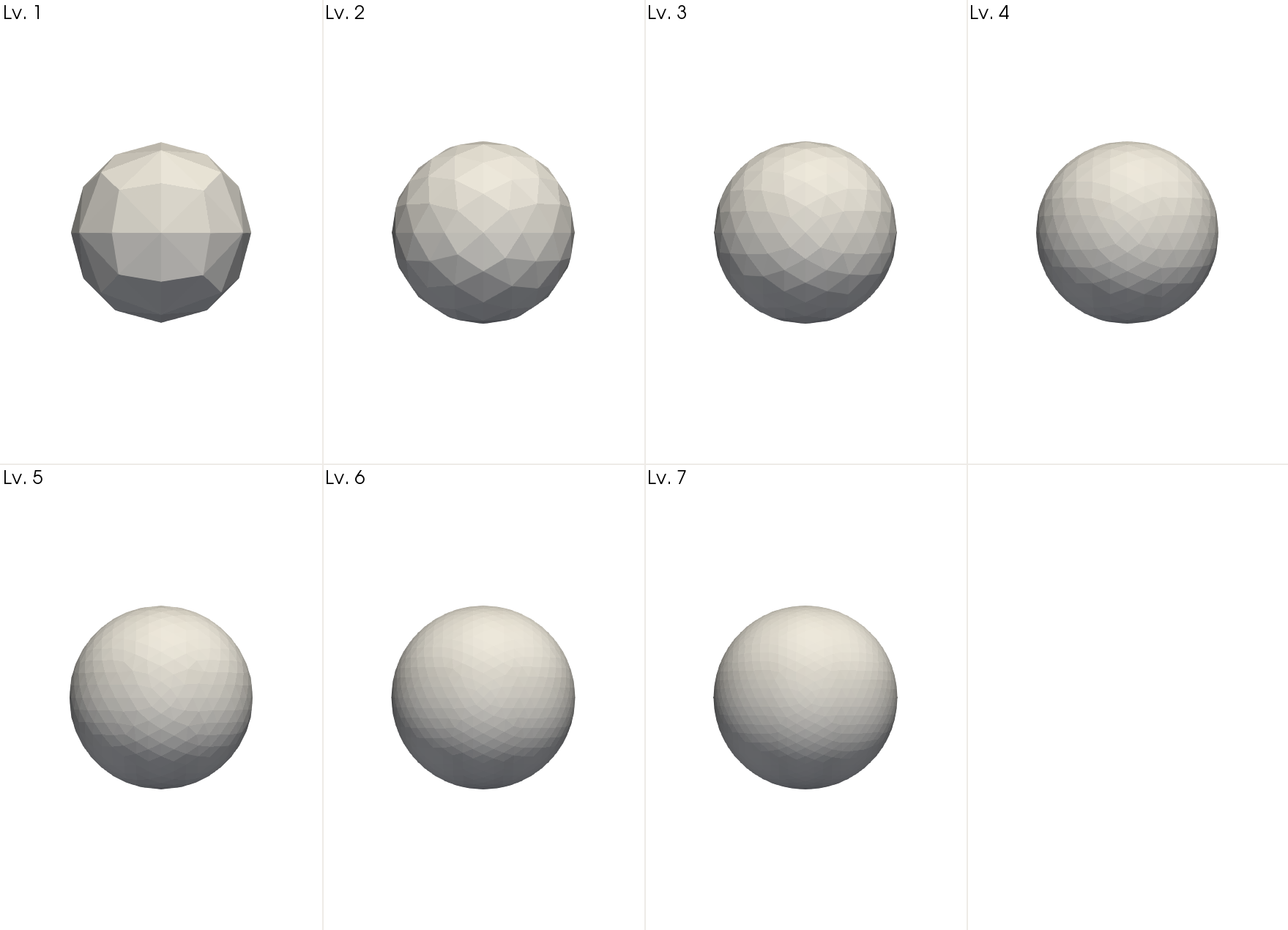}}
  \caption{Mesh refinement levels \label{fig:sphere_collage_with_text}} 
\end{subfigure}
\caption{Spheres (a) arranged in Fibonacci spiral configuration and (b) showing the different mesh refinement levels used.  \label{fig:sphere_problem_setup}}
\end{figure}

The physical setup of the problem consists of a central, heated sphere surrounded by several ambient temperature spheres of the same size arranged in a spiral pattern, as shown in Fig.~\ref{fig:sphere_problem_setup}. Referring to the labeling system shown in Fig.~\ref{fig:spiral_overhead_view_labels}, sphere 1 is located at the origin, and then spheres 2, 3, 5, 7, 9, 11, and 13 are offset by $(\pm a_n, \pm a_n)$, from the previous sphere's location, where $a_n$ is the $n$th Fibonacci number, and the signs are chosen so that a counter-clockwise spiral is formed. Finally, the remaining five spheres are placed at the center of the quarter-circular arcs joining the previous spheres. The spheres are meshed with four-node tetrahedral elements, and have 88, 240, 432, 720, 1248, 1992, and 2664 exterior facets at Levels 1--7, respectively, as shown in Fig.~\ref{fig:sphere_collage_with_text}. The total number of mesh nodes, elements, and cavity facets (including all 13 spheres) for each mesh level is given in Table~\ref{tab:mesh_levels}.

This geometric configuration provides several different opportunities for low-rank approximation. First, because the distance between the spheres increases rapidly as one proceeds around the spiral, there will be cavity facets in relatively close proximity (near the center of the spiral) as well as far apart, relative to the mesh spacing (facet sizes) present in the mesh. Second, because the spheres are arranged in a planar spiral, there are some facets (around 7\% of the total number of facets) on the ``top'' and ``bottom'' of each sphere, as well as along the exterior of the outer arm of the spiral, which cannot ``see'' any other facets, i.e.\ are isolated.

\begin{table}[htb]
    \caption{Number of facets, nodes, and elements for each mesh level in the Fibonacci spheres problem. We use a first-order Lagrange finite element, so the number of degrees of freedom (DOFs) in the temperature solve is equal to the number of nodes for each mesh level.\label{tab:mesh_levels}}
    \centering
    \begin{tabular}{cccc}
    \toprule
    Mesh Level & Num.\ Facets & Num.\ Nodes & Num. Elements \\
    \midrule
         1 & \num{1144 } & \num{793  } & \num{2782  } \\
         2 & \num{3120 } & \num{3081 } & \num{14339 } \\
         3 & \num{5616 } & \num{7553 } & \num{38766 } \\
         4 & \num{9360 } & \num{13910} & \num{72956 } \\
         5 & \num{16224} & \num{32097} & \num{175422} \\
         6 & \num{25896} & \num{59241} & \num{328211} \\
         7 & \num{34632} & \num{90051} & \num{505258} \\
    \bottomrule
    \end{tabular}
\end{table}

Although this configuration of spheres does not conform to the classical definition of a ``cavity,'' i.e.\ a convex, closed region, the characteristics listed previously make it an interesting test problem for assessing the effectiveness of the low-rank solver's capabilities. The difference between a closed and open cavity is reflected in the treatment of the view factor matrix entries as follows. First, we define the (scaled) $i$th row sum
\begin{align}
    \label{eq:rowsum}
    s_i := \frac{1}{A_i} \sum_j F_{ij}
\end{align}
For the ``truth'' view factor matrix, we of course have $0 \leq s_i \leq 1$, but because of quadrature error and floating point tolerances, we sometimes obtain row sums slightly outside this range. Therefore, we define the corresponding ``clamped'' row sum
\begin{align}
    \label{eq:clamped_rowsum}
    c_i := \left\{
\begin{array}{cc}
     0, &  s_i < 0 \\
     1, & s_i > 1 \\
     s_i, & \text{otherwise}\\
\end{array}
    \right.
\end{align}
Then, in the case of a closed cavity, the $i$th row of the view factor matrix is scaled by $c_i$, i.e.\ if $c_i \neq 0$,
\begin{align}
    \label{eqn:F_closed}
    F_{ij} \leftarrow \frac{F_{ij}}{c_i} \quad \forall i,j
\end{align}
In an open cavity, on the other hand, we do not scale the view factor matrix rows, but instead define an ``ambient'' temperature, $T_{\infty}$, and include an additional nonlinear boundary flux term in the residual \eqref{eq:SVF_Galerkin} which is given by
\begin{align}
    \label{eqn:flux_to_ambient}
    B_N(T_h) := \sum_i \int_{\Gamma_i} \left(1 - c_i\right) \left(T_h^4 - T_{\infty}^4\right) \varphi_N \, \dx{S}
\end{align}
to represent thermal energy which is ``radiated to ambient.''

While the most ``realistic'' approach would be to treat this model as an open cavity problem, with almost all the thermal energy being radiated to ambient, and little heat actually being transferred between the spheres, we found that by instead treating this model as an artificially (mathematically) closed cavity using the approach described above, the solutions had a more pronounced transient behavior and hence provided a better test case for our implementation. That is, in the open cavity case, the surrounding spheres' temperature remains nearly unchanged, while in the closed cavity case, the surrounding spheres are heated non-uniformly. This choice allowed us to better classify the accuracy of the different low-rank approaches relative to the Direct solve approach.

The initial temperature of all the spheres is set to 300K, except for the central sphere, which is set to 1000K. There is no internal heat generation, and the radiation cavity is defined to be the union of the surfaces of all the spheres. The final solution time is set to 1000s, and we use a uniform time step of $\Delta t = 25$s, so that a total of 40 time steps are computed. We use the same time discretization for all mesh levels, though as mentioned elsewhere, a more realistic mesh refinement study would refine in both time and space simultaneously to reduce both temporal and spatial discretization errors at approximately the same rate. We treat this model as a closed cavity, which means that all thermal energy radiated from the central sphere is eventually absorbed (and possibly re-emitted) by the other spheres.

\begin{figure}[htb]
\centering
\includegraphics[width=0.9\textwidth]{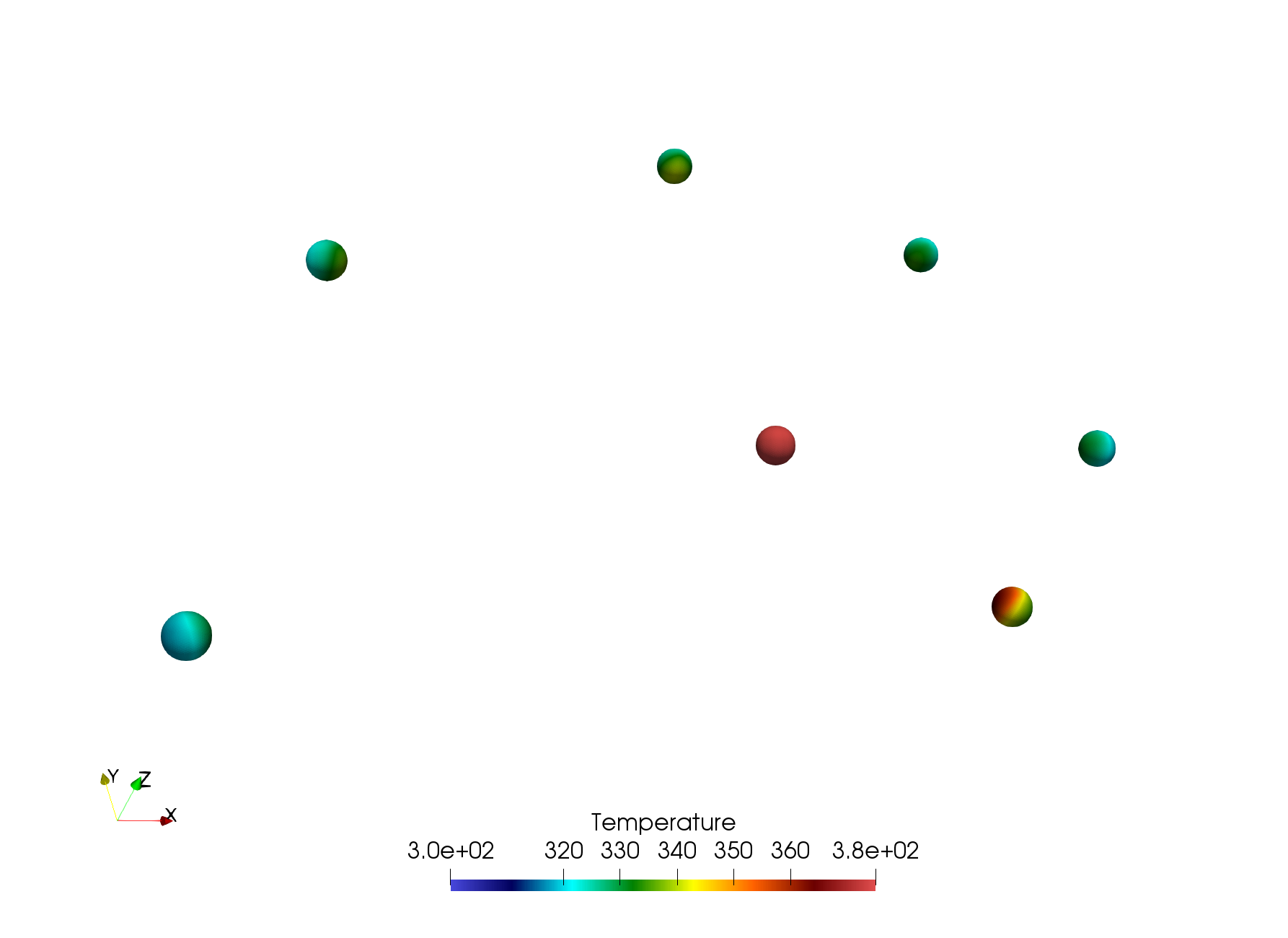}
\caption{Surface temperature of central spheres at time $t=1000s$.\label{fig:close_in_spheres}}
\end{figure}

At steady state, all spheres should reach the same equilibrium temperature. A representative solution at time $t=1000s$ on the Level 7 mesh is shown in Fig.~\ref{fig:close_in_spheres}, where the color bar is set to a maximum of 380K in order to better highlight the temperature variation on the non-central spheres; this causes the central sphere to appear as a single color which is outside the maximum temperature range. Also, note that we are only showing the central seven spheres in Fig.~\ref{fig:close_in_spheres} since the outer spheres are both too far away to show in the same figure with a reasonable level of detail, and the temperature fields on the outer spheres are more uniform and hence less interesting than the inner spheres.


\begin{figure}[htb]
\centering
\begin{subfigure}[t]{.5\textwidth}
  \centering
  \includegraphics[width=.95\linewidth]{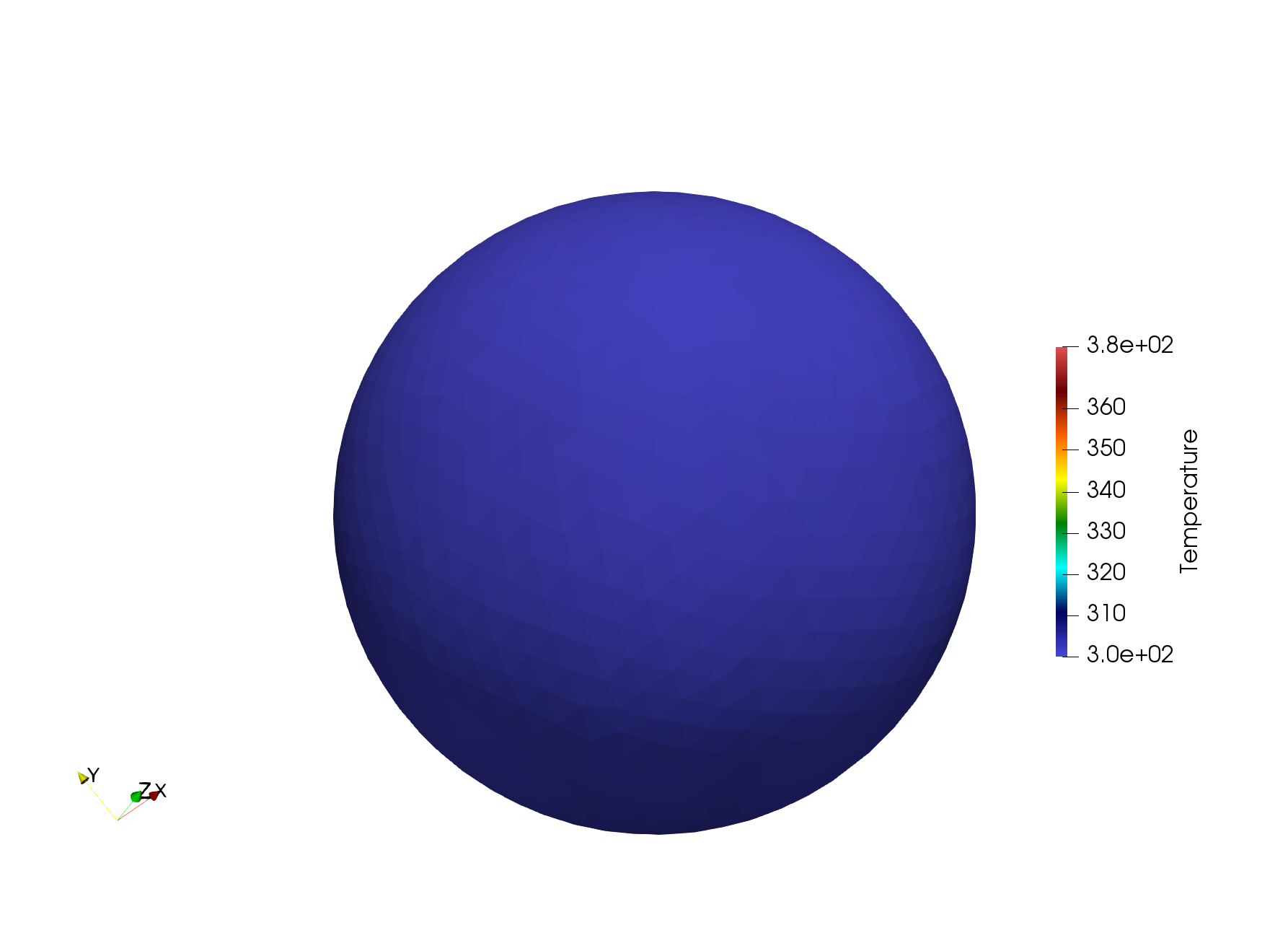}
  \caption{$t=0$s \label{fig:close_sphere_a}}
 \end{subfigure}%
\begin{subfigure}[t]{.5\textwidth}
  \centering
  \includegraphics[width=.95\linewidth]{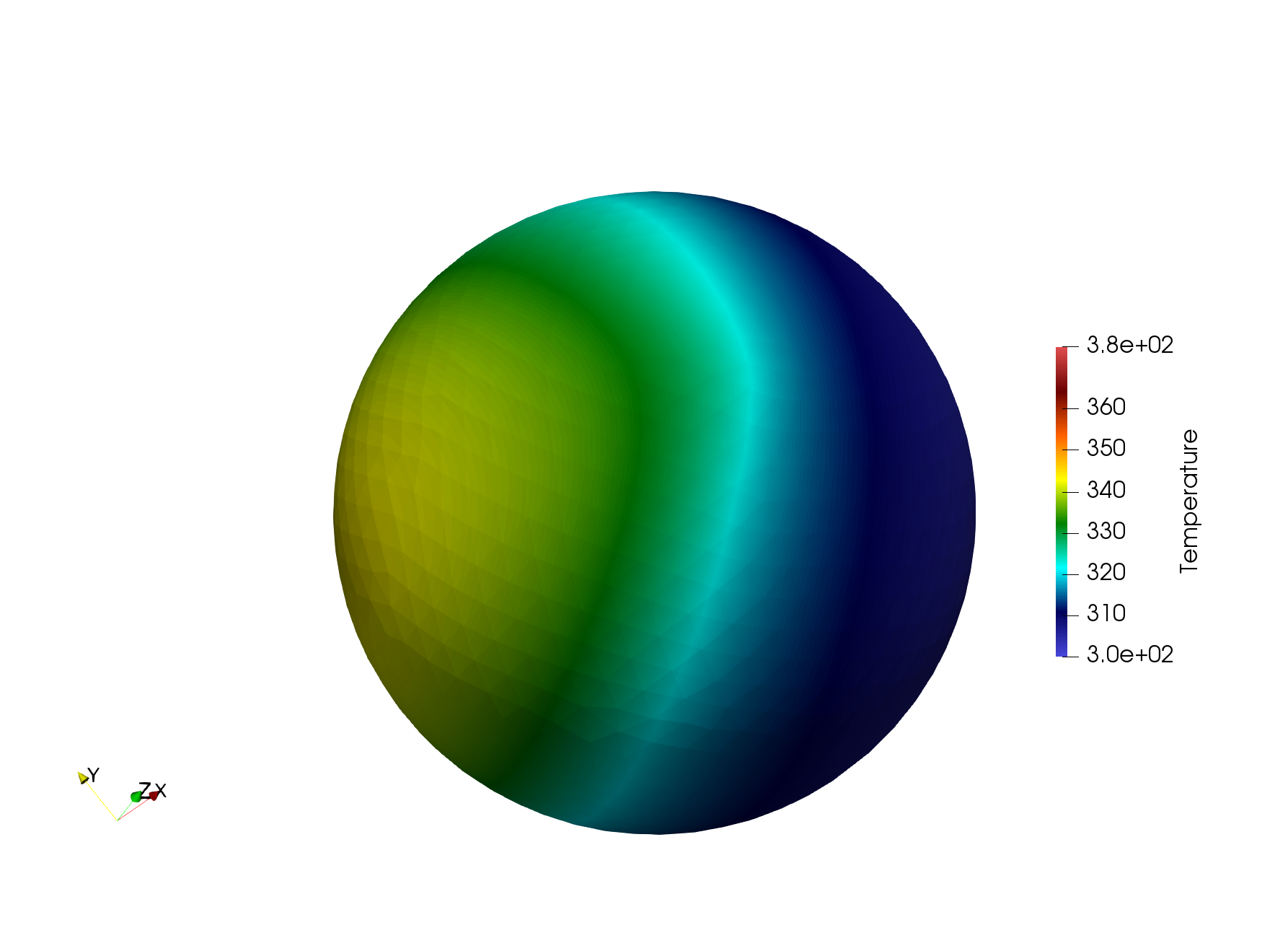}
  \caption{$t=200$s \label{fig:close_sphere_b}} 
\end{subfigure}
\\
\begin{subfigure}[t]{.5\textwidth}
  \centering
  \includegraphics[width=.95\linewidth]{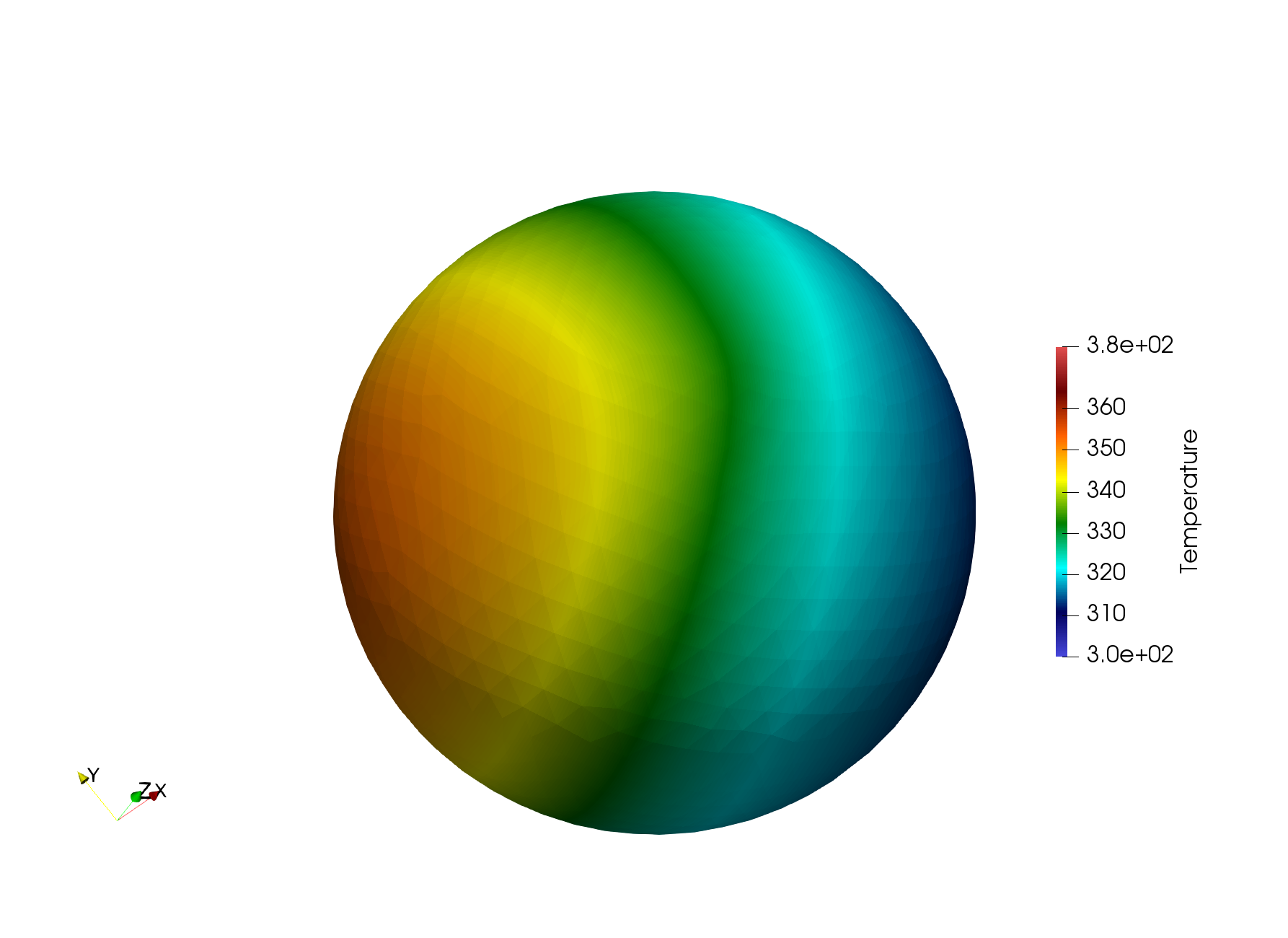}
  \caption{$t=400$s \label{fig:close_sphere_c}} 
 \end{subfigure}%
\begin{subfigure}[t]{.5\textwidth}
  \centering
  \includegraphics[width=.95\linewidth]{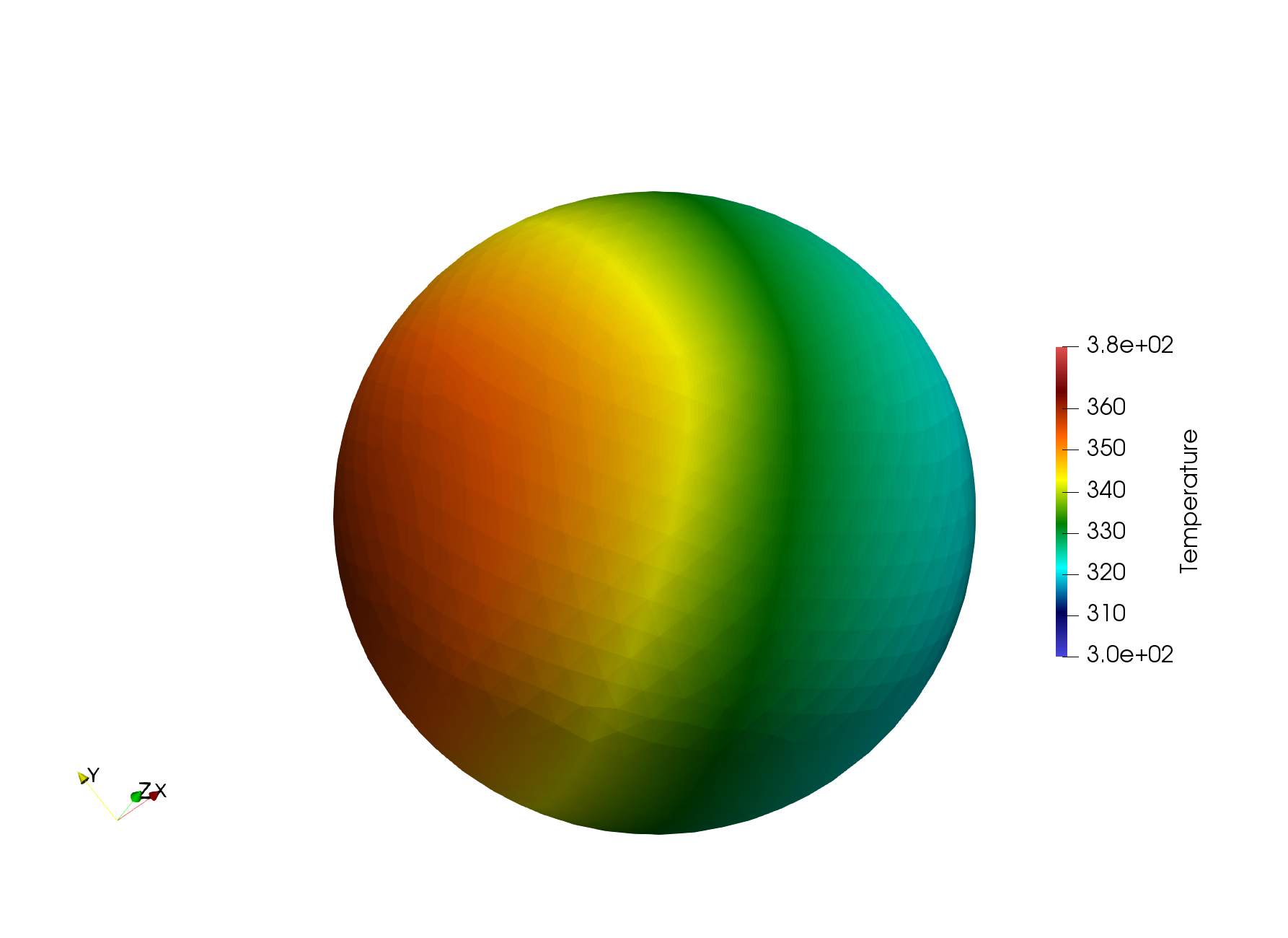}
  \caption{$t = 600$s \label{fig:close_sphere_d}} 
\end{subfigure}
\\
\begin{subfigure}[t]{.5\textwidth}
  \centering
  \includegraphics[width=.95\linewidth]{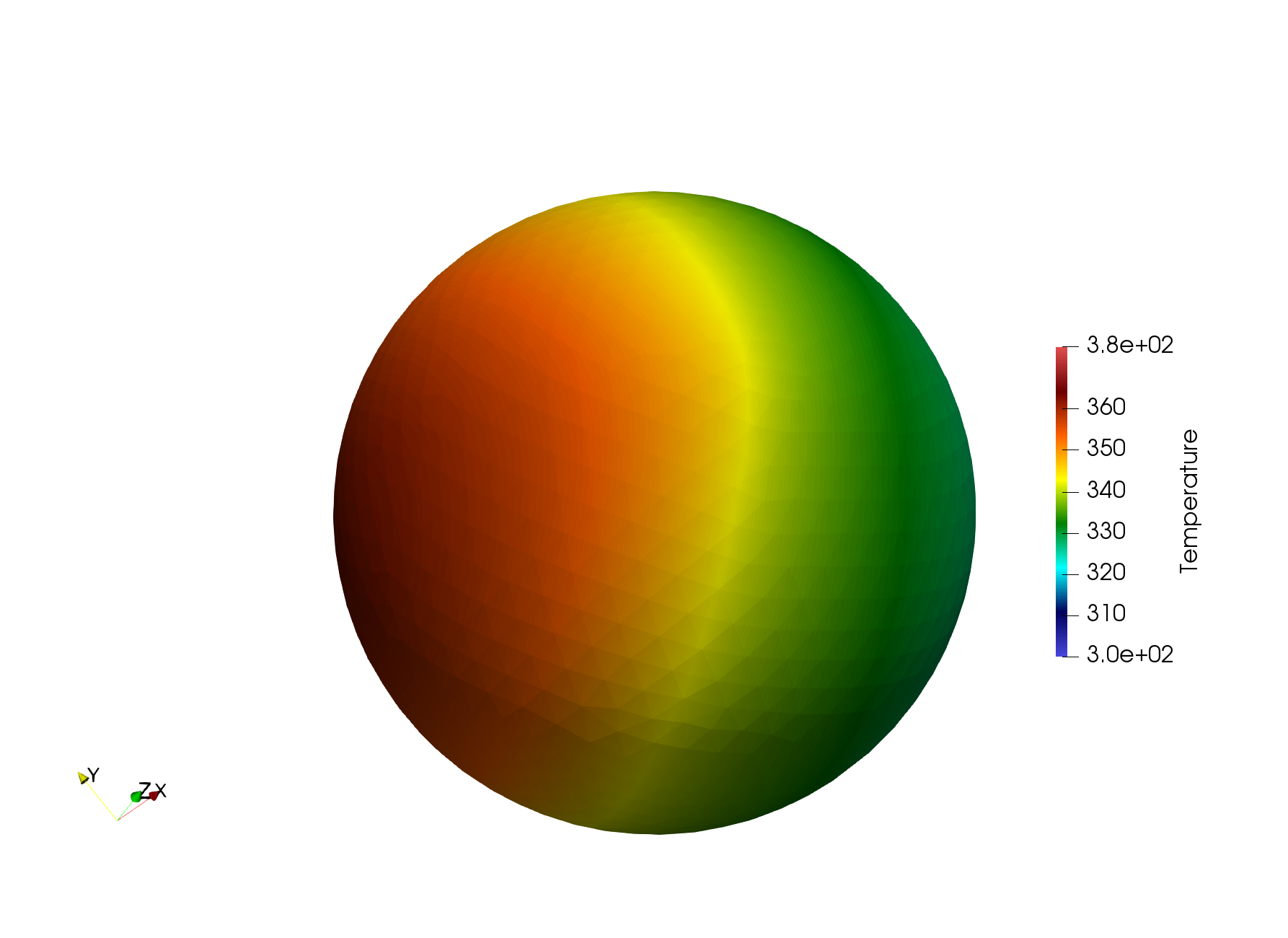}
  \caption{$t = 800$s \label{fig:close_sphere_e}} 
 \end{subfigure}%
\begin{subfigure}[t]{.5\textwidth}
  \centering
  \includegraphics[width=.95\linewidth]{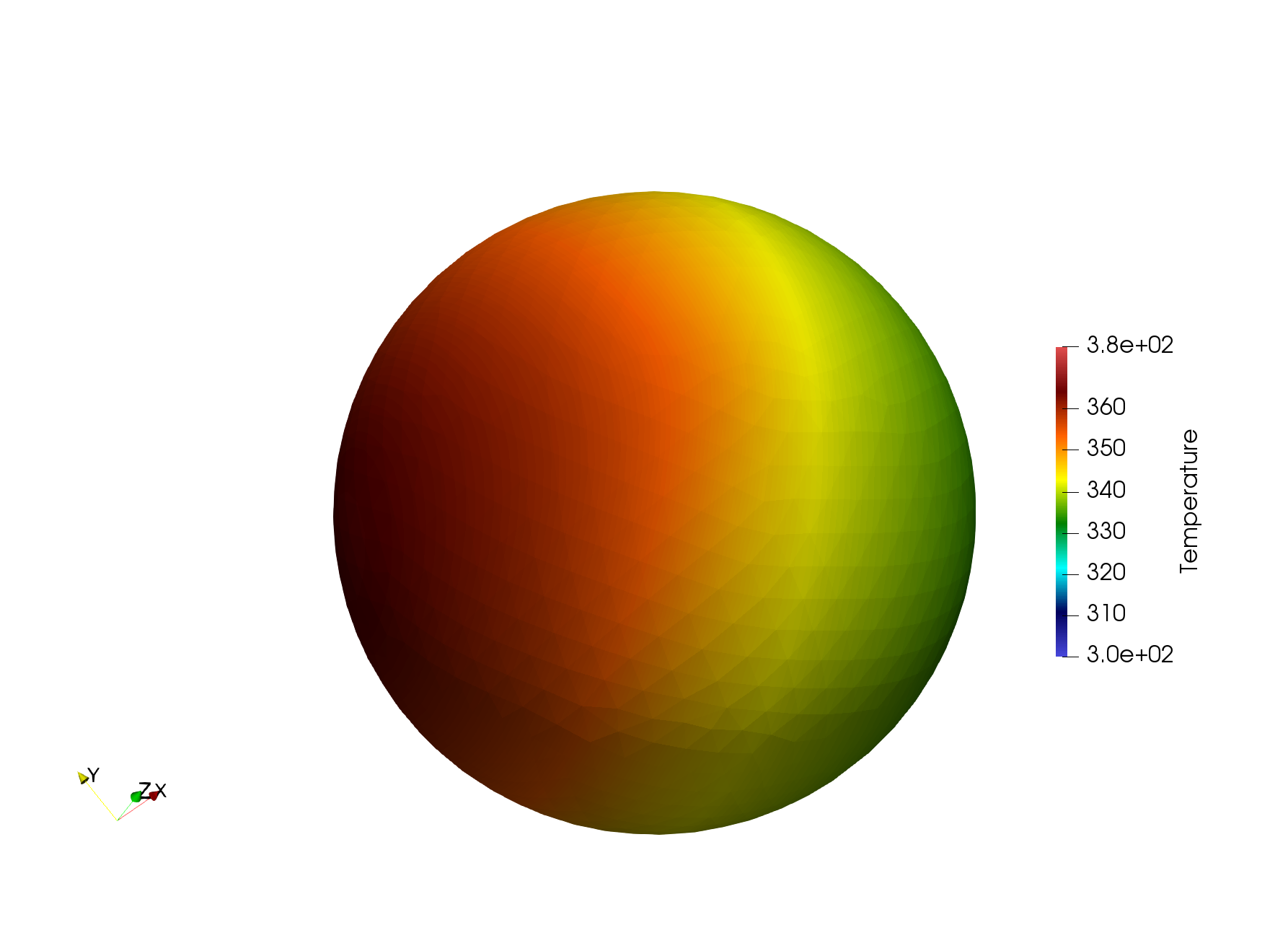}
  \caption{$t = 1000$s \label{fig:close_sphere_f}} 
\end{subfigure}
\caption{Surface temperature of sphere 2 for various simulation times.\label{fig:close_sphere_temps}}
\end{figure}

\begin{figure}[htb]
\centering
\begin{subfigure}[t]{.5\textwidth}
  \centering
  \includegraphics[width=.95\linewidth]{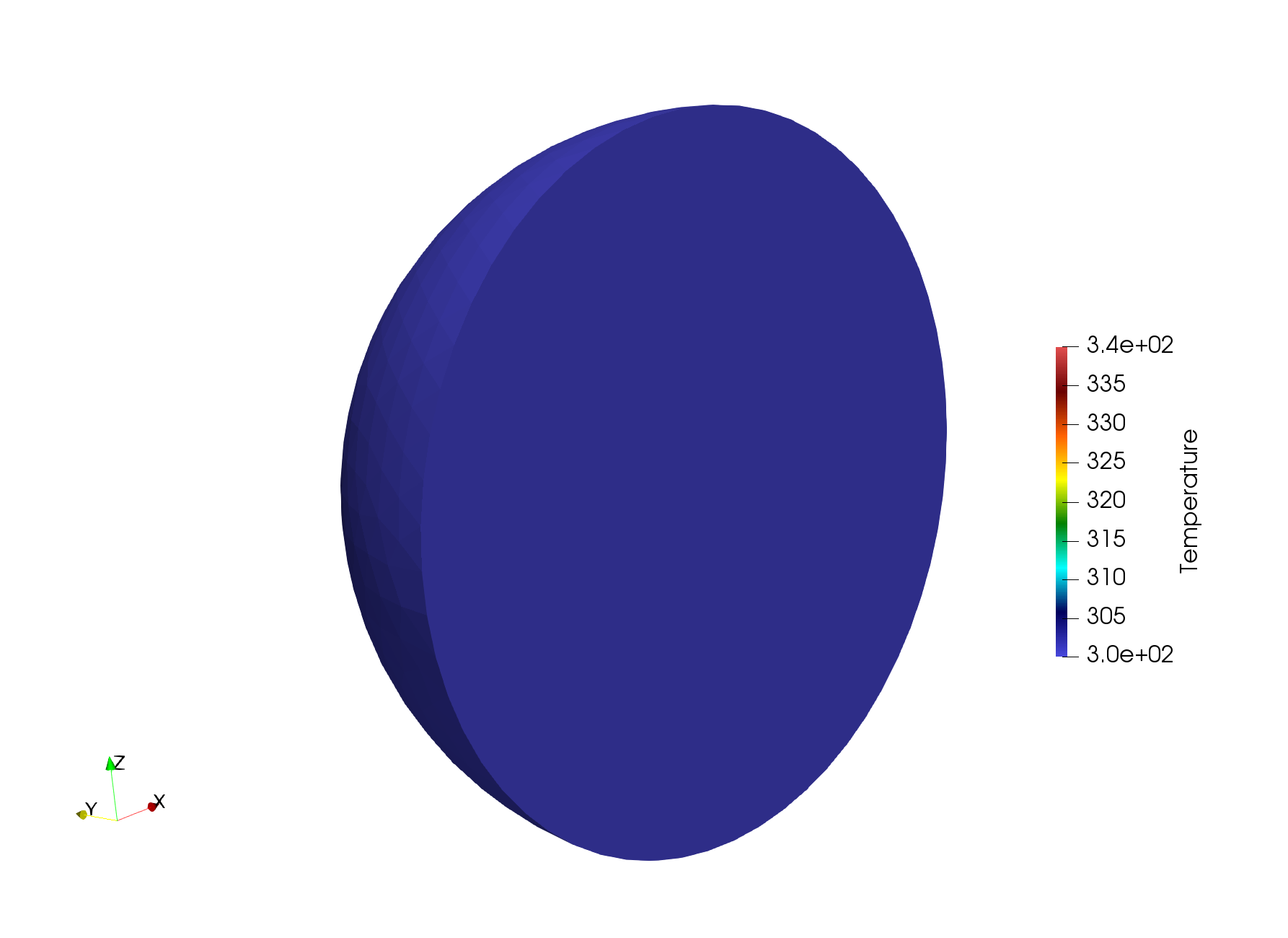}
  \caption{$t=0$s \label{fig:sphere_3_a}}
 \end{subfigure}%
\begin{subfigure}[t]{.5\textwidth}
  \centering
  \includegraphics[width=.95\linewidth]{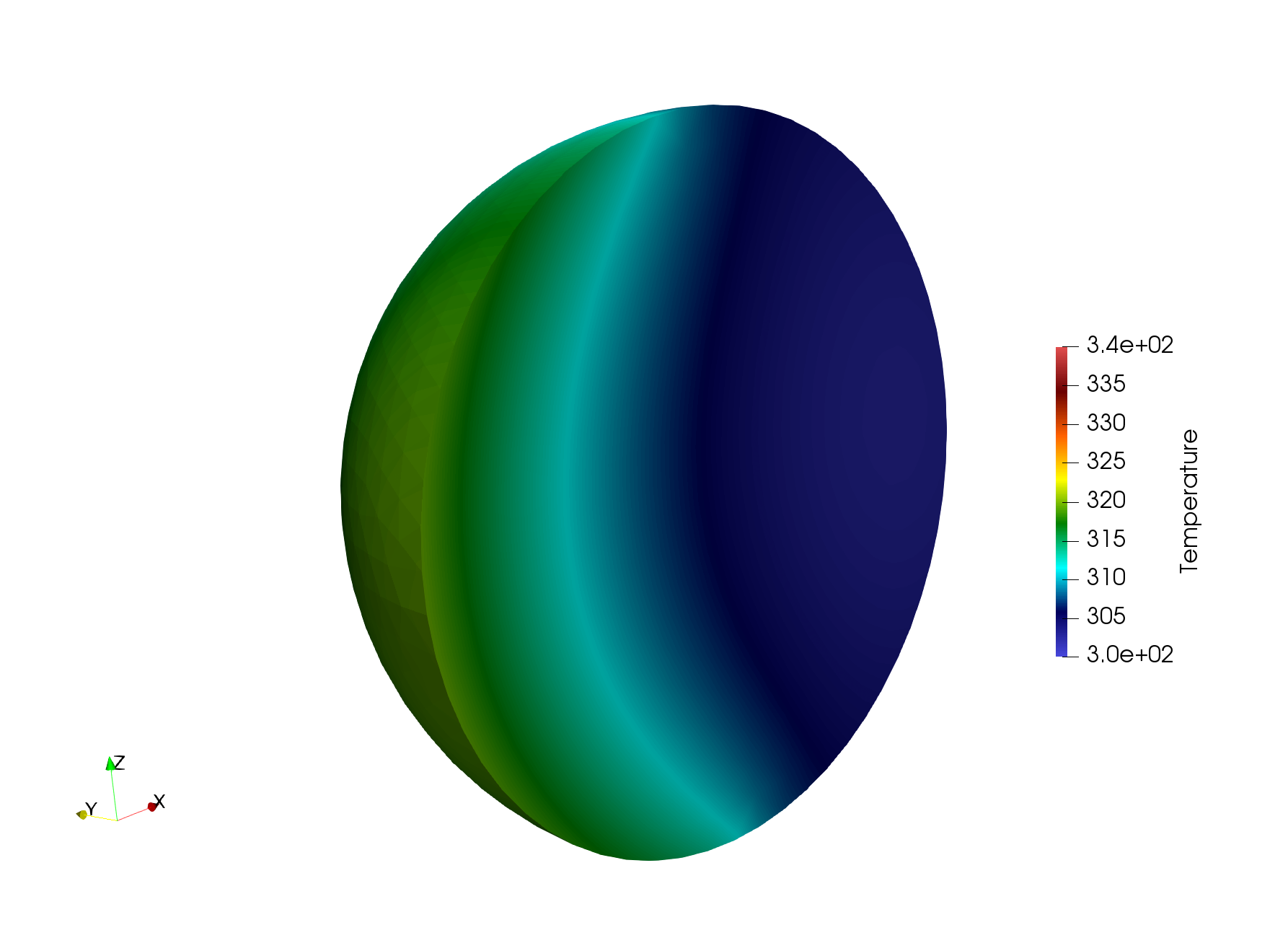}
  \caption{$t=200$s \label{fig:sphere_3_b}} 
\end{subfigure}
\\
\begin{subfigure}[t]{.5\textwidth}
  \centering
  \includegraphics[width=.95\linewidth]{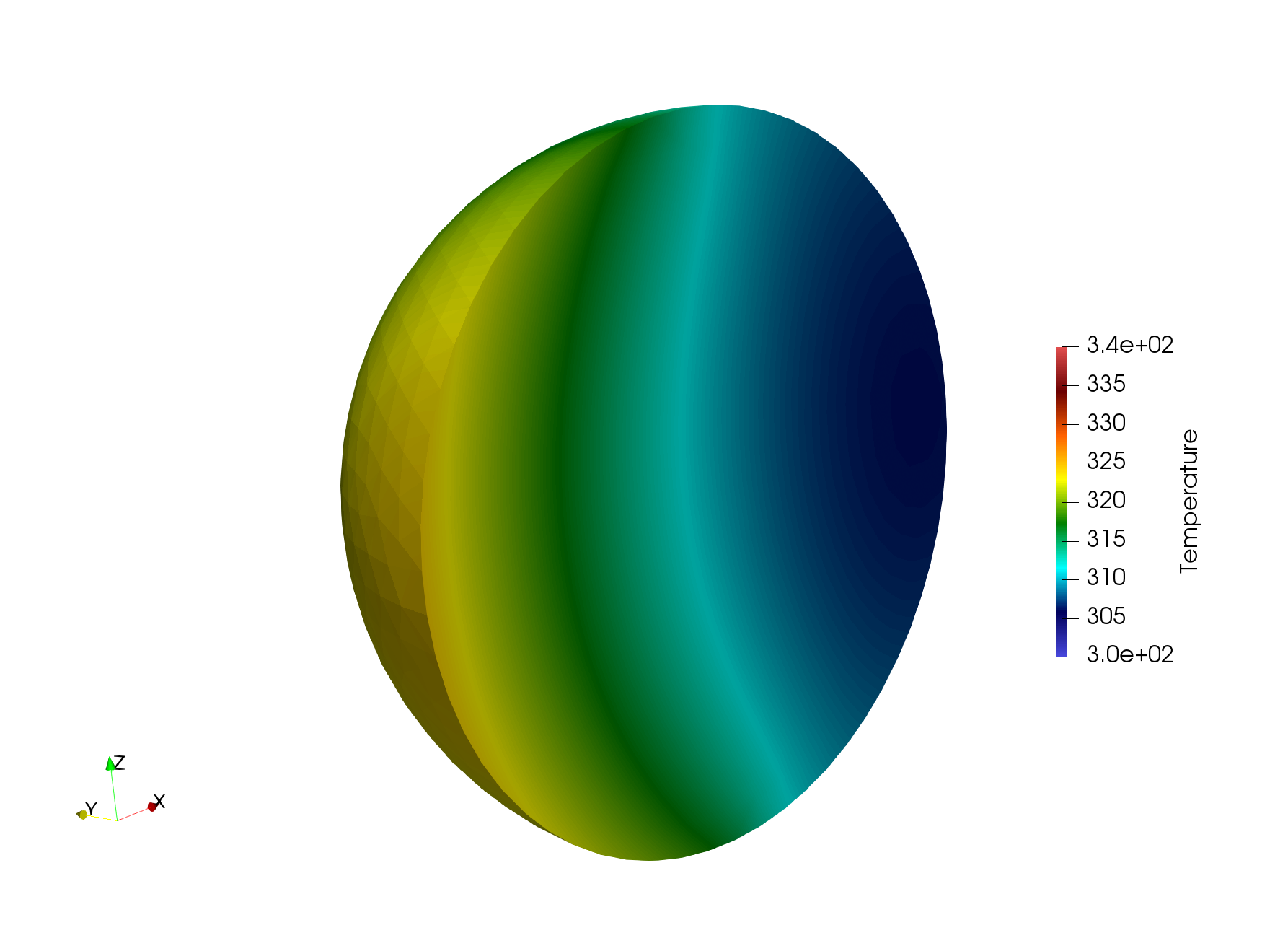}
  \caption{$t=400$s \label{fig:sphere_3_c}} 
 \end{subfigure}%
\begin{subfigure}[t]{.5\textwidth}
  \centering
  \includegraphics[width=.95\linewidth]{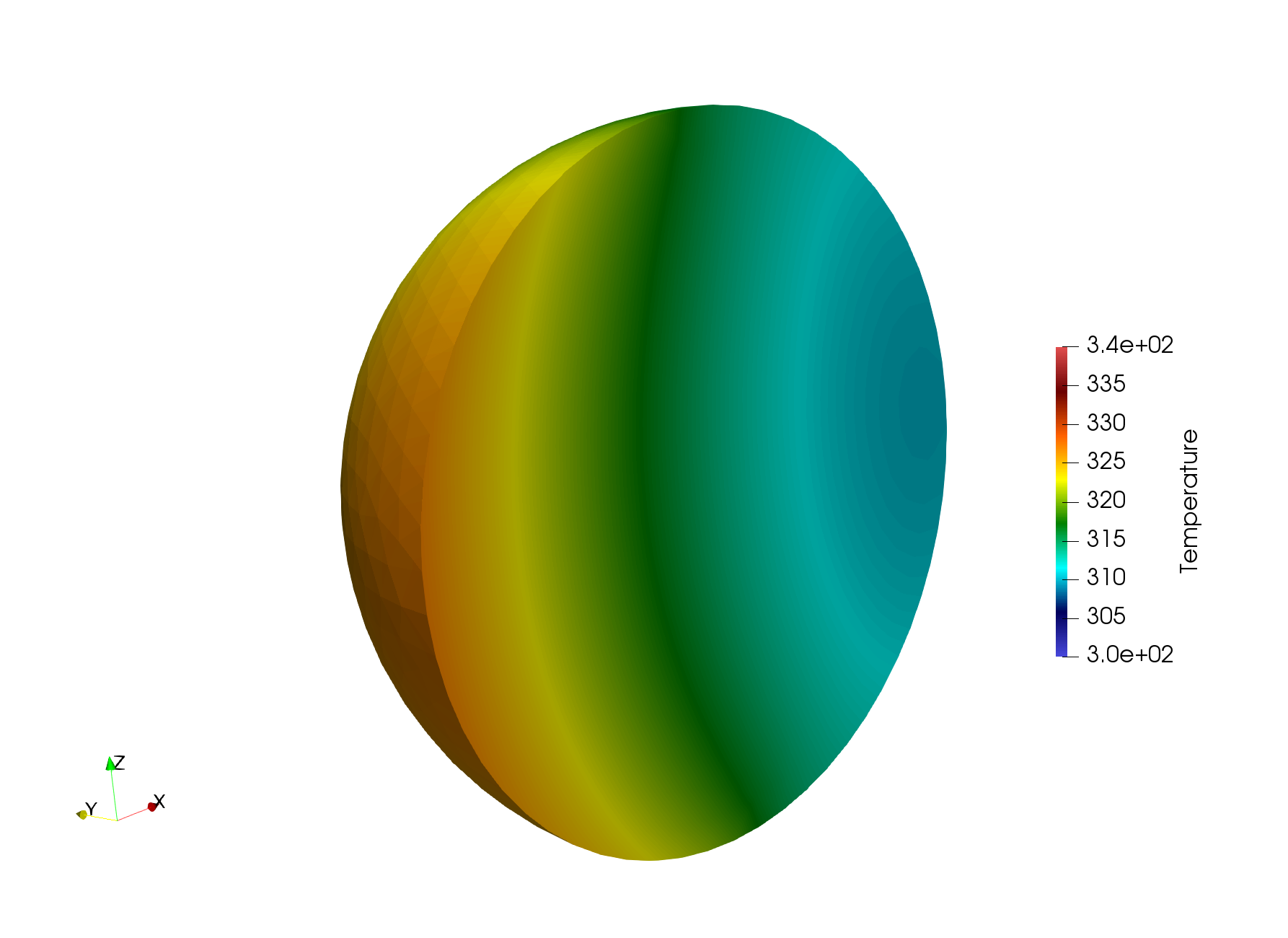}
  \caption{$t = 600$s \label{fig:sphere_3_d}} 
\end{subfigure}
\\
\begin{subfigure}[t]{.5\textwidth}
  \centering
  \includegraphics[width=.95\linewidth]{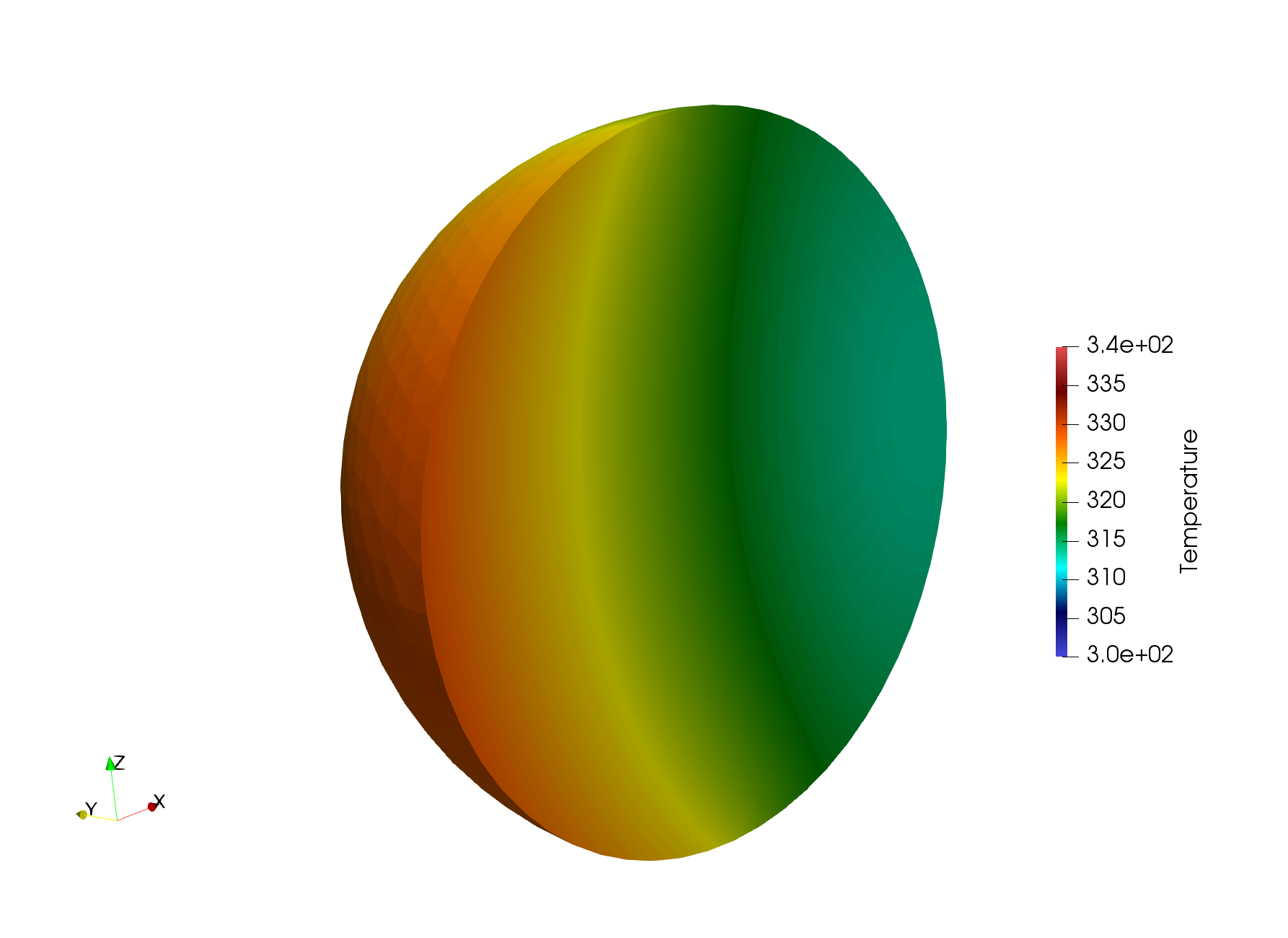}
  \caption{$t = 800$s \label{fig:sphere_3_e}} 
 \end{subfigure}%
\begin{subfigure}[t]{.5\textwidth}
  \centering
  \includegraphics[width=.95\linewidth]{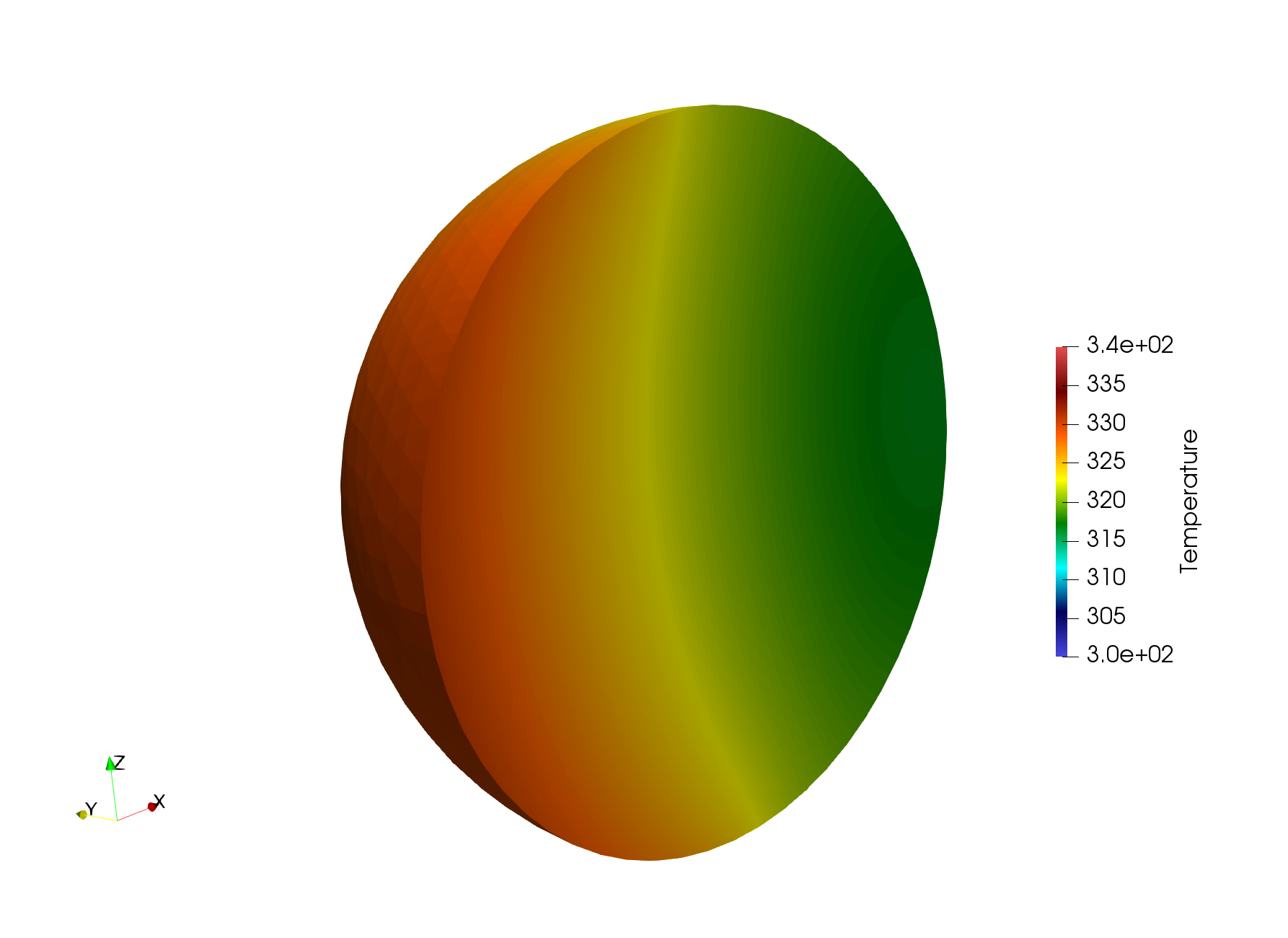}
  \caption{$t = 1000$s \label{fig:sphere_3_f}} 
\end{subfigure}
\caption{Clip plane showing internal temperature of sphere 3 for various simulation times.\label{fig:sphere_3_temps}}
\end{figure}

More detailed images of sphere 2 (the sphere closest to the central sphere) at evenly-spaced time steps are shown in Fig.~\ref{fig:close_sphere_temps}. This figure clearly shows that the side of sphere 2 which is closest to the central sphere is preferentially heated relative to the far side, which remains several degrees cooler. In Fig.~\ref{fig:sphere_3_temps}, we show the internal details of the temperature field in sphere 3, by creating a clip plane whose normal is orthogonal to the line segment which joins the centers of spheres 1 and 3. This figure shows how heat is transferred internally in the sphere (via conduction) after the surface is heated via radiation. The temperature field contours within the sphere are curved rather than straight, a result of the non-uniform radiation heat flux being applied to the surface.

\clearpage
\subsection{Accuracy comparison with Direct solve method\label{sec:comparison}}
When possible, we compare the performance of the low-rank solver to the so-called ``Direct solve'' method in which the full view factor and reflection matrices are explicitly constructed, and the reflection matrix inverse is explicitly computed.  That is, the $\Jcav$ contribution from \eqref{eq:Jproduct} is a true dense matrix, and the Direct solver thus has no way of taking advantage of e.g.\ the zero rows arising from the isolated facets mentioned previously.

In our comparisons below in which the low-rank truncation tolerance $\epsrel$ is varied, the Direct solution method is taken as the ``truth'' solution, and the relative error in the low-rank solution methods is computed via the discrete space-time norm
\begin{align}
    \label{eq:eT}
    e(T) := \max_{t_n} \frac{\| T - T_{\text{Direct}}\|}{\|T_{\text{Direct}}\|}
\end{align}
where $T_{\text{Direct}}$ is the solution computed by the Direct solve method, the maximum is taken over all time steps $t_n$, $n=1,\ldots,n_{\text{steps}}$, and $\|\cdot \|$ is the discrete $\ell^2$-norm. In addition to comparing the temperature fields computed by the Direct and low-rank approaches, we also provide detailed comparisons of the overall runtime and memory usage of the different approaches, in order to quantify the effectiveness of the low-rank solver approach.

\begin{figure}[!htb]
\centering
\begin{subfigure}[t]{.5\textwidth}
  \centering
  \includegraphics[width=\linewidth]{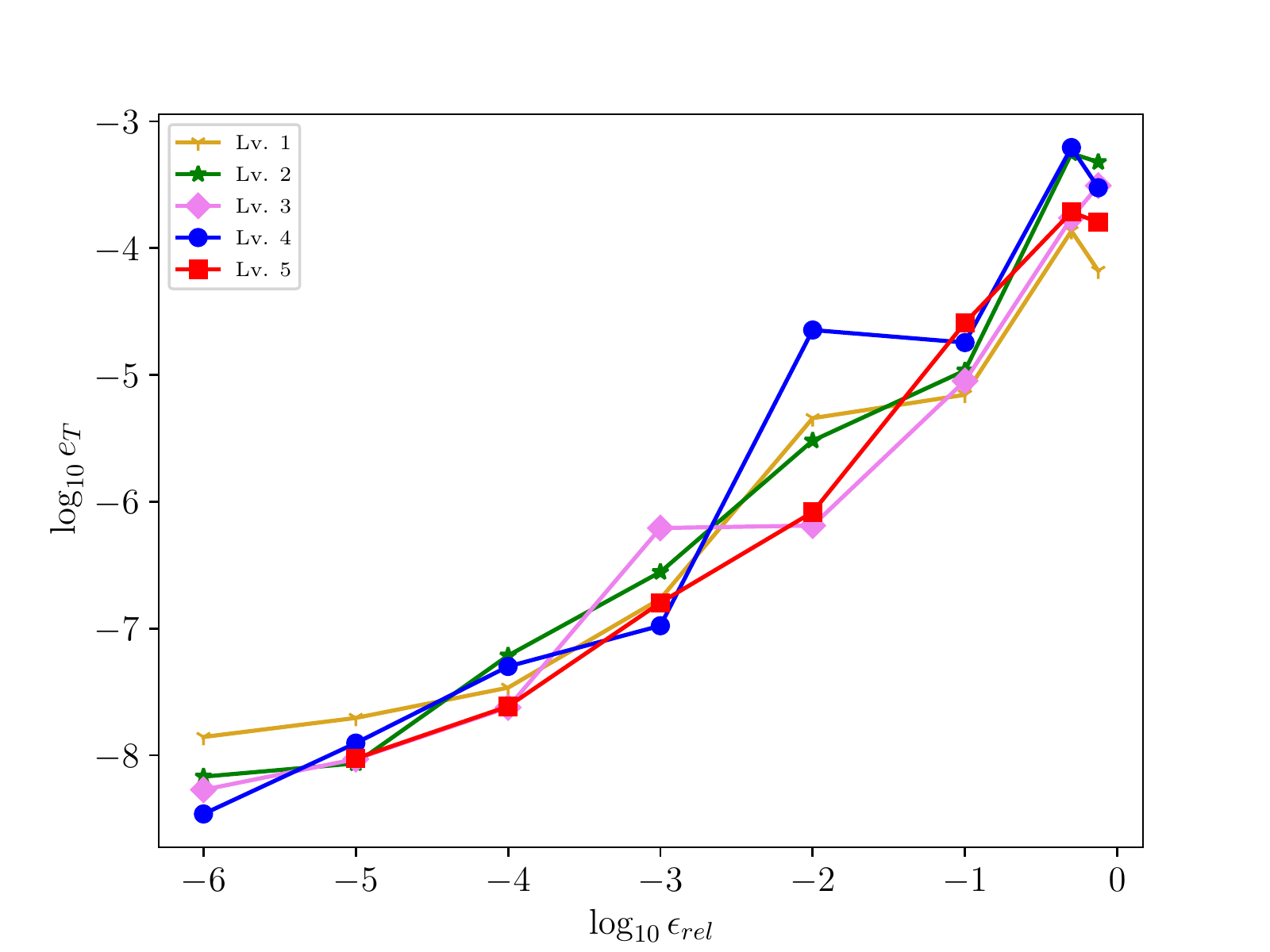}
  \caption{Relative error $e(T)$ \label{fig:fib_13_spheres_full_piv_ACA_rel_err}}
 \end{subfigure}%
\begin{subfigure}[t]{.5\textwidth}
  \centering
  \includegraphics[width=\linewidth]{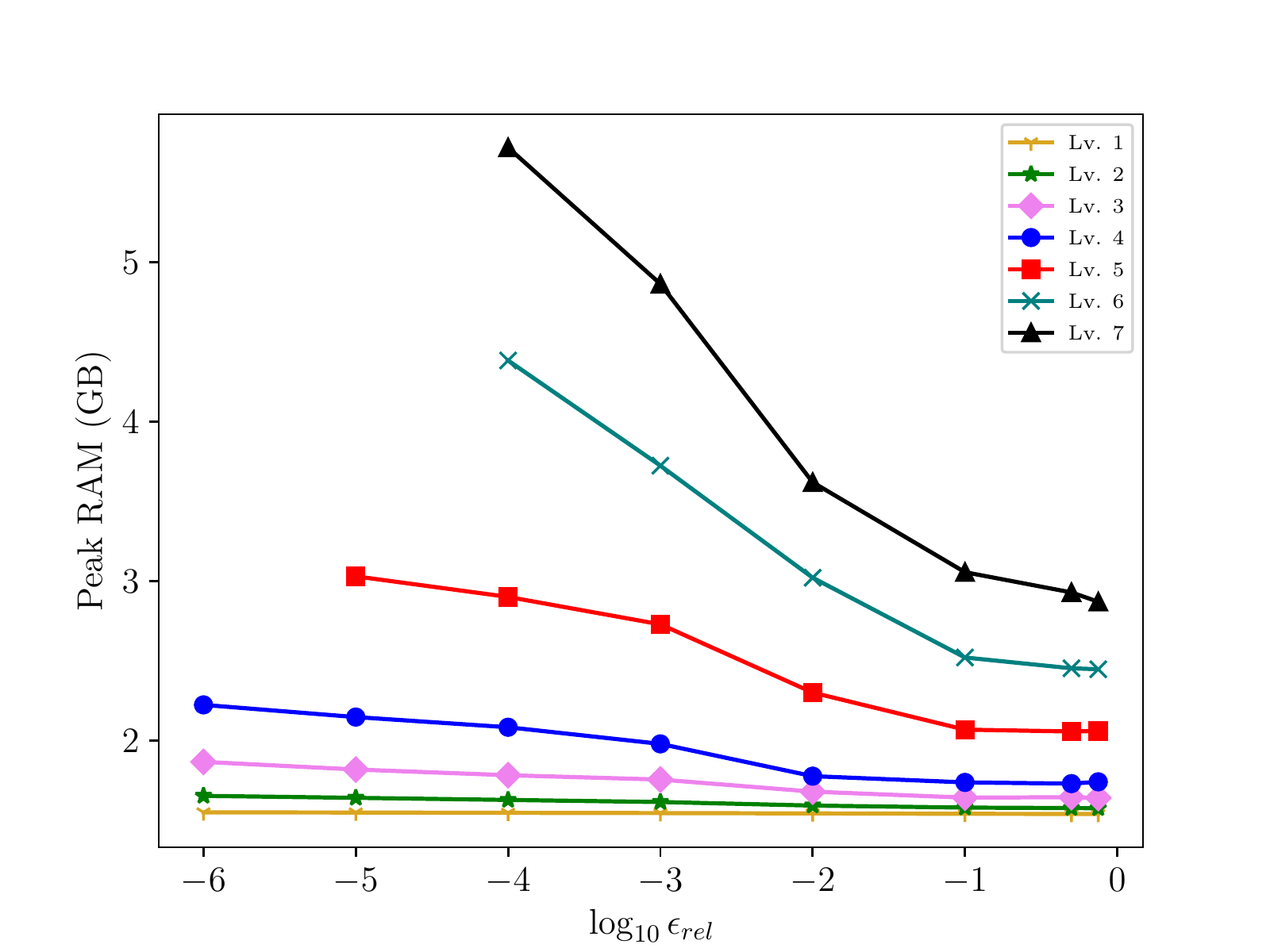}
  \caption{Peak RAM usage\label{fig:fib_13_spheres_peak_ram}} 
\end{subfigure}
\caption{Plots comparing (a) relative error $e(T)$ defined in \eqref{eq:eT} and (b) peak RAM usage vs. $\epsrel$ for all mesh refinement levels.\label{fig:fib_13_spheres_full_piv_ACA_speedup_and_peak_ram}}
\end{figure}

In our comparisons we employed the sequence of meshes from Level 1 to 5 shown in Fig.~\ref{fig:sphere_collage_with_text}. On each mesh considered, we also varied the truncation tolerance $\epsrel$ (cf. \eqref{eq:erel}) from a maximum value of 0.75 to a minimum value of $10^{-6}$. The results of these calculations are shown in Fig.~\ref{fig:fib_13_spheres_full_piv_ACA_rel_err}, and provide good evidence that the accuracy of the low-rank solver depends continuously on the truncation accuracy $\epsrel$, but does not depend strongly on the mesh refinement level. In practice, this means that one can choose a single value of $\epsrel$ for an entire mesh refinement study, and expect a consistent level of agreement with the corresponding Direct solve at that mesh refinement level.

The relative error comparisons stop at Level 5 since the Direct solver required more RAM than was available on the compute node in order to successfully solve the problem on the Level 6 and 7 meshes. The low-rank solver, on the other hand, had relatively modest peak memory requirements, as shown in Fig.~\ref{fig:fib_13_spheres_peak_ram}. For the low-rank solver, the peak memory usage increases quickly for decreasing $\epsrel$ on fine grids, but for the largest $\epsrel$ values (coarsest truncation accuracy), we observe that the peak RAM usage for the Level 7 mesh is only about 1.87 times larger than for the Level 1 mesh, despite the Level 7 mesh having over 30 times as many facets in the cavity.

\clearpage
\subsection{The choice of $\epsrel$\label{sec:epsrel_choice}}
Due to the importance of the choice of $\epsrel$, we make some further comments about this tolerance here. Recall that this tolerance refers to the truncation accuracy of the data-sparse blocks within the block low-rank approximation, whereas the dense blocks are represented exactly. This means that in general we expect a higher-accuracy approximation of $F$ and $C$ than the level specified by $\epsrel$, given that we have a mix of ($\epsrel$ approximate) low-rank and (exact) dense blocks. We quantify this effect in Table~\ref{tab:F_err}, which shows that the approximation error in $F$ is typically an order of magnitude lower than $\epsrel$.

\begin{table}[htb]
    \caption{Relative error in the low-rank approximation, $F_{\epsrel}$, with respect to the ``true'' view factor matrix, $F$, for different $\epsrel$ values. $\|\cdot\|$ represents the Frobenius norm.\label{tab:F_err}}
    \centering
    \begin{tabular}{ll}
    \toprule
    $\epsrel$ & $\|F - F_{\epsrel}\| / \|F\|$ \\
    \midrule
    \num{1.e-1} &   \num{9.529e-03}   \\ 
    \num{1.e-2} &   \num{9.668e-04} \\ 
    \num{1.e-3} &   \num{9.416e-05} \\
    \num{1.e-4} &   \num{9.399e-06} \\
    \num{1.e-5} &   \num{9.372e-07} \\
    \num{1.e-6} &   \num{9.415e-08}  \\
    \bottomrule
    \end{tabular}
\end{table}

A second point we note is that we are ultimately only interested in the accuracy of the solution field, $T$, and the cavity radiation terms are simply a ``stepping stone'' to $T$. This is an important point because the solution fields in heat transfer are well known to be smooth, which tends to compensate for approximation introduced in the representation of $F$ and $C$. Taken together, these two effects explain why in general we observe $e(T) \ll \epsrel$, and hence this provides justification for choosing relatively large values of $\epsrel$ such as $\epsrel=0.1$.

This is welcome since these relatively large values of $\epsrel$ are in the region in which we observe good scalability and high speed from the block low-rank approach in our numerical results. Nevertheless, it is important to keep in mind that in general the behavior of the block low-rank approximation with respect to truncation tolerance will be problem-dependent, and hence one should ideally justify a particular choice of $\epsrel$ by first running some test cases targeted at the model problem of interest.

\clearpage
\subsection{Performance vs.\ $\epsrel$\label{sec:comparison_epsrel}}

The performance of the low-rank solver can roughly be categorized into two main parts: (i) the time spent in computing and applying the low-rank LU factorization, and (ii) everything else including finite element assembly, GMRES iterations, file I/O, and other tasks. Here we will mainly be interested in the time spent in part (i), since part (ii) represents tasks that must also be performed in the Direct solver. Within part (i), it is useful to further separate the time spent (a) building the view factor matrix, (b) computing the low-rank LU factorization, and (c) applying the low-rank LU factorization. Items (a) and (b) happen only once per simulation, and thus their costs are amortized across all time steps, while the time required for item (c) increases in direct proportion to the number of time steps computed.

\begin{figure}[!htb]
\centering
\begin{subfigure}[t]{.5\textwidth}
  \centering
  \includegraphics[width=\linewidth]{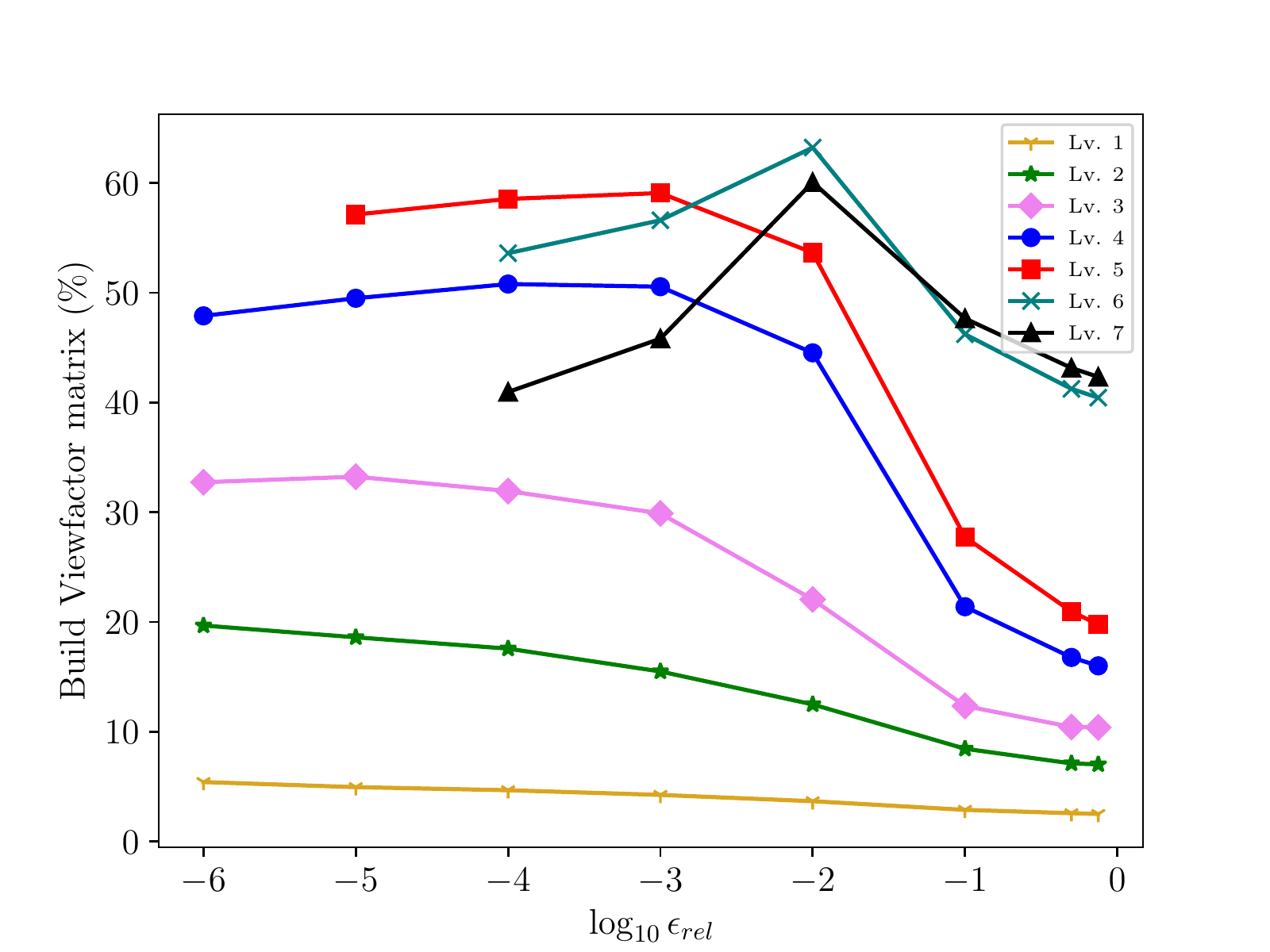}
  \caption{Build view factor matrix \label{fig:fib_13_spheres_full_piv_ACA_build_vf_time}}
 \end{subfigure}%
\begin{subfigure}[t]{.5\textwidth}
  \centering
  \includegraphics[width=\linewidth]{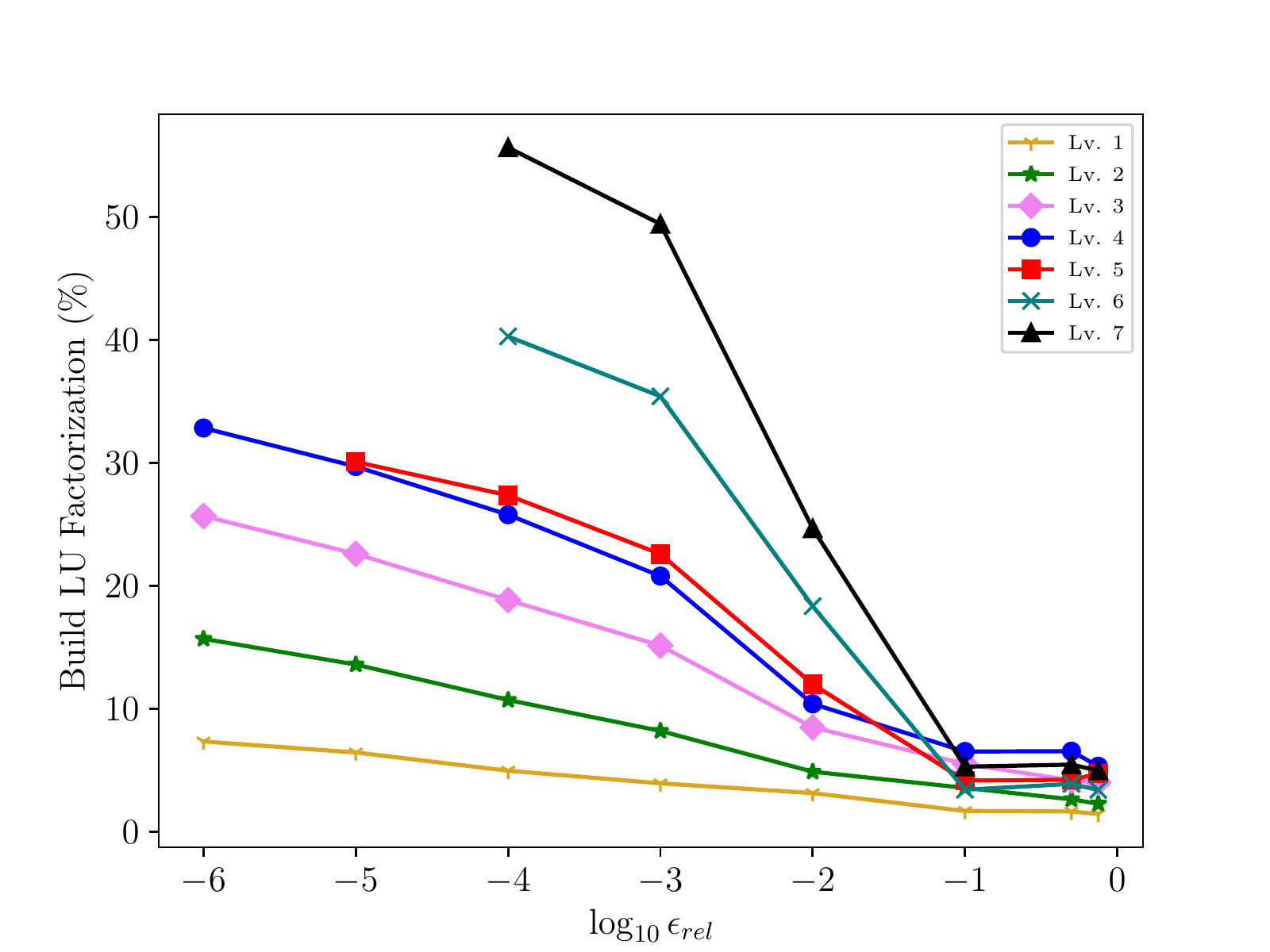}
  \caption{Build LU factorization \label{fig:fib_13_spheres_full_piv_ACA_LU_time}}
\end{subfigure}
\\
\begin{subfigure}[t]{.5\textwidth}
  \centering
  \includegraphics[width=\linewidth]{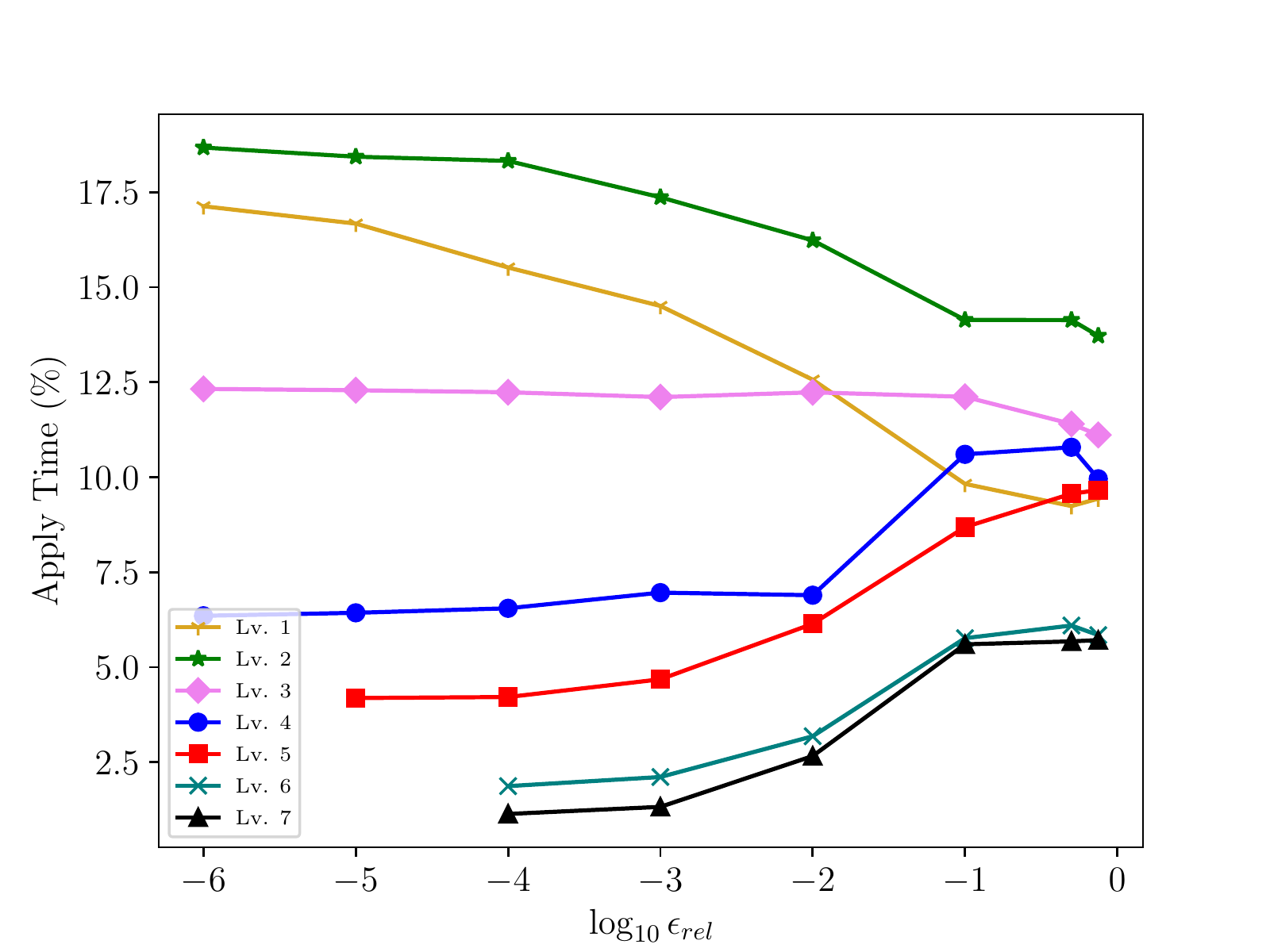}
  \caption{Apply LU factorization \label{fig:fib_13_spheres_full_piv_ACA_apply_time}}
\end{subfigure}
\caption{Time required to (a) build the view factor matrix, (b) build the low-rank LU factorization of the cavity Jacobian, and (c) apply the low-rank LU factorization (as a percentage of the total runtime) vs.\ $\epsrel$ for each mesh level.\label{fig:fib_13_spheres_full_piv_ACA_factorize_and_apply}}
\end{figure}

The time spent in part (i), i.e.\ building the view factor matrix, computing the low-rank factorization, and applying the LU factorization is shown, for a range of $\epsrel$ values, in Fig.~\ref{fig:fib_13_spheres_full_piv_ACA_factorize_and_apply}. For coarse meshes, we find that as the truncation tolerance $\epsrel$ is decreased, the percentage of time spent building the view factor matrix (Fig.~\ref{fig:fib_13_spheres_full_piv_ACA_build_vf_time}) plateaus or even decreases, while the percentage of time spent computing the LU factorization (Fig.~\ref{fig:fib_13_spheres_full_piv_ACA_LU_time}) increases steadily. In general, decreasing the truncation tolerance $\epsrel$ has a more pronounced effect on the LU factorization time for the finer (Level 6 and 7) grids; the LU factorization time exceeds 50\% of the total runtime in the finest case.

Our test problem only computes 40 time steps, so we would of course expect the percentage of time spent building the view factor matrix and LU factorization to be lower in a model that performed hundreds or even thousands of time steps. Regarding the time spent applying the LU factorization, for the coarse mesh levels tested, the percentage of time spent applying the factorization increases with decreasing $\epsrel$, but on finer meshes it actually decreases with decreasing $\epsrel$. This indicates the relatively larger expense of building the low-rank LU factorization vs.\ applying it for fine grids. In problems with more time steps, we would expect the application time of the low-rank factorization to eventually dominate the creation time.

\clearpage
\subsection{Performance vs.\ Direct solver\label{sec:performance_direct}}

\begin{figure}[htb]
\centering
\begin{subfigure}[t]{.5\textwidth}
  \centering
  \includegraphics[width=\linewidth]{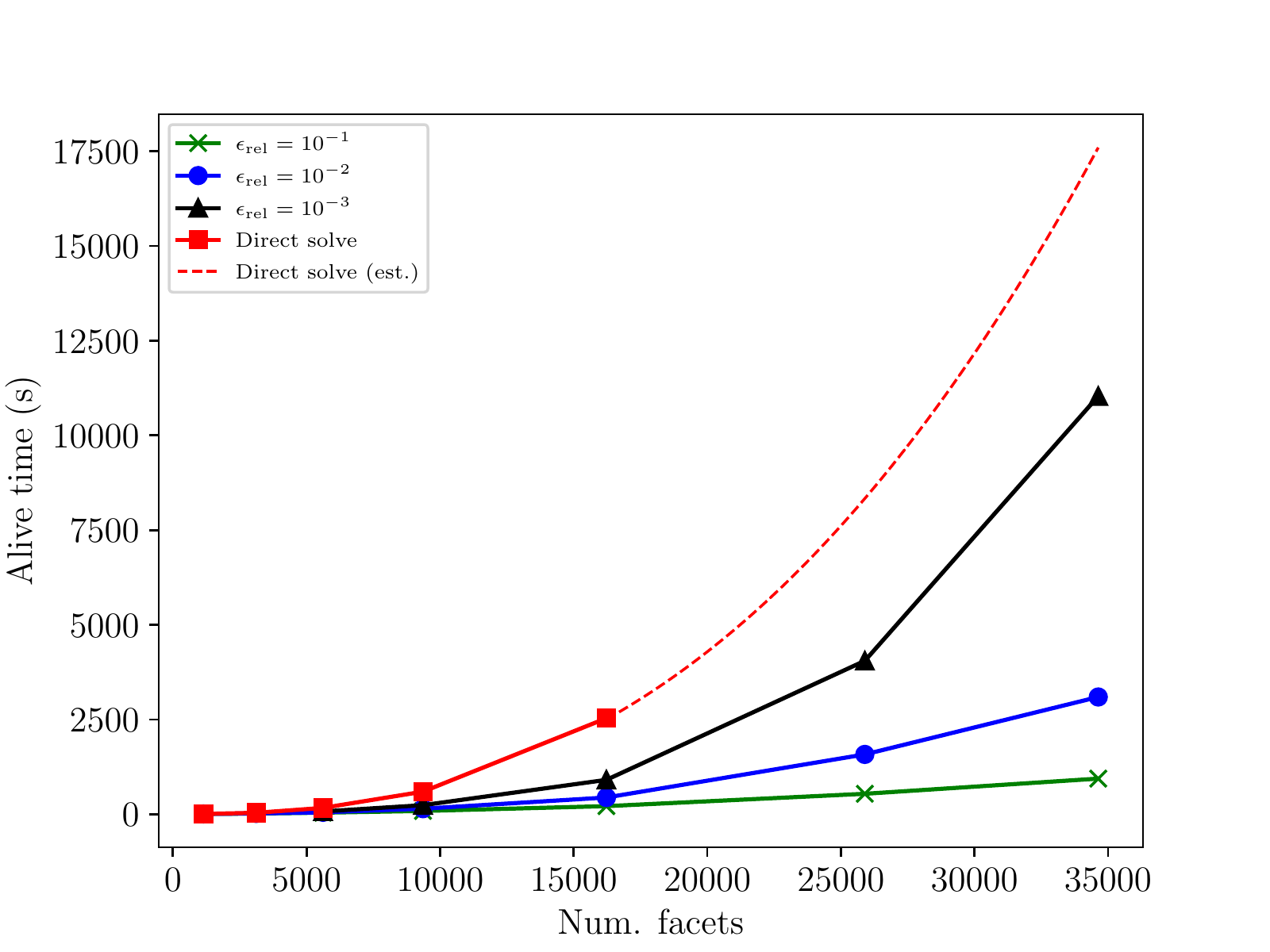}
  \caption{Alive time \label{fig:fib_13_spheres_full_piv_ACA_alive_times}}
 \end{subfigure}%
\begin{subfigure}[t]{.5\textwidth}
  \centering
  \includegraphics[width=\linewidth]{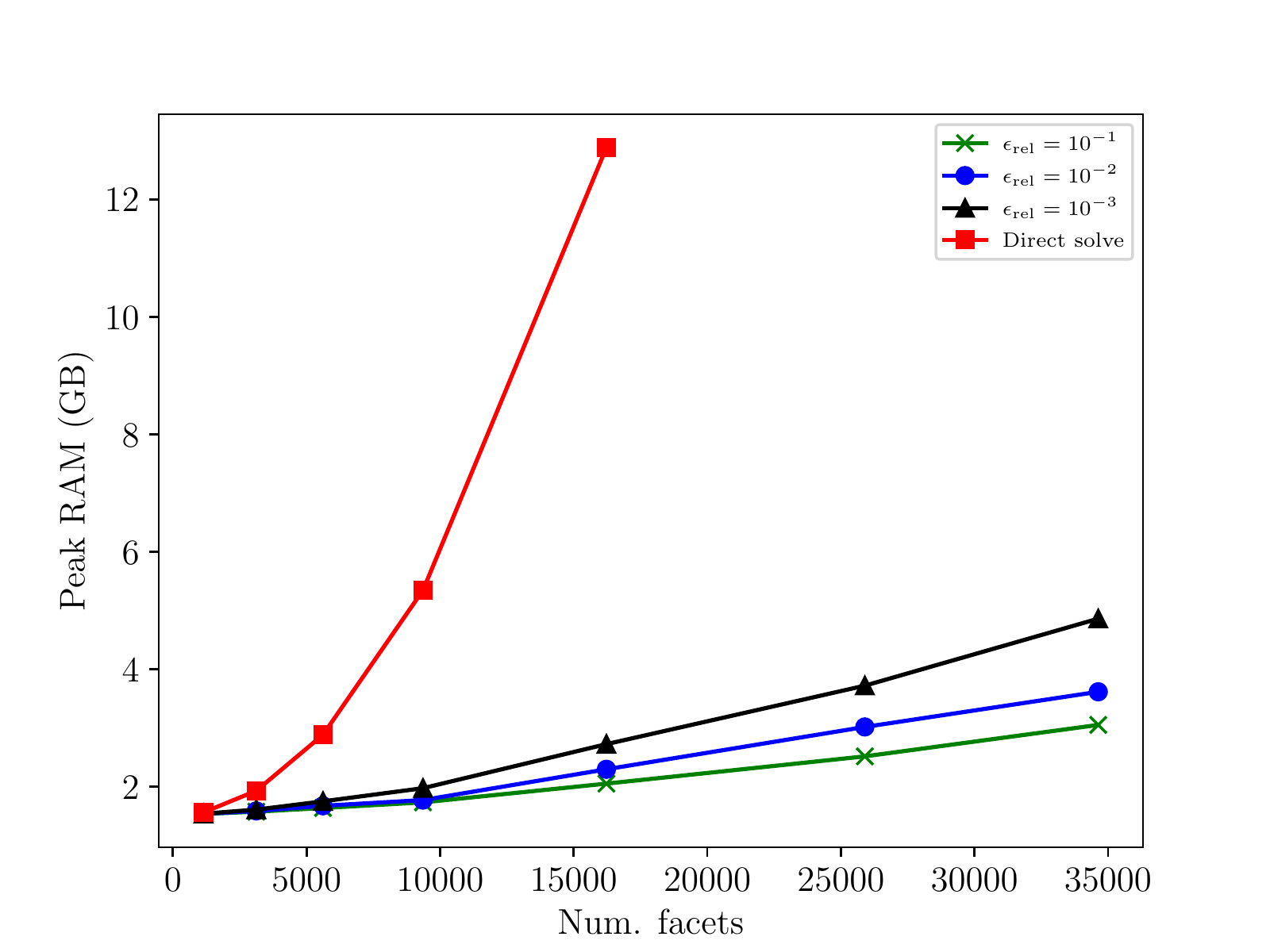}
  \caption{Peak RAM usage\label{fig:fib_13_spheres_full_piv_ACA_peak_ram}} 
\end{subfigure}
\caption{Comparison of (a) total simulation time and (b) peak RAM usage for the Direct and low-rank solvers with $\epsrel= 10^{-1}, 10^{-2}, 10^{-3}$ vs.\ number of radiation cavity facets.\label{fig:fib_13_spheres_full_piv_ACA_alive_time_and_ram}}
\end{figure}

The previous plots have focused on comparing the low-rank solver's performance for different mesh refinement levels and different $\epsrel$ values, but it is of course also necessary to compare the performance of the low-rank solver to the Direct solve method. In Fig.~\ref{fig:fib_13_spheres_full_piv_ACA_alive_time_and_ram}, we therefore switch the $x$-axis of the plots from $\epsrel$ to the number of cavity facets, a quantity which is the same for both the Direct and low-rank solver at all mesh levels. In Fig.~\ref{fig:fib_13_spheres_full_piv_ACA_alive_times}, we compare the full simulation time (also referred to as ``Alive Time'') for the Direct solver and the low-rank solver with three different $\epsrel$ values of $10^{-1}$, $10^{-2}$, and $10^{-3}$.

As mentioned previously, we were unable to run the Direct solver on the Level 6 and 7 meshes due to the RAM required to form the explicit cavity Jacobian inverse, and this rapid memory growth is depicted in Fig.~\ref{fig:fib_13_spheres_full_piv_ACA_peak_ram}. Moreover, the dashed line in Fig.~\ref{fig:fib_13_spheres_full_piv_ACA_alive_time_and_ram} indicates extrapolated Direct solve simulation times for the Level 6 and 7 meshes, based on a curve fit of the preceding data points. For the Level 7 mesh with $\epsrel = 10^{-1}$, the low-rank solver's simulation time is approximately 15.6 minutes while the estimated time for the direct solver (on a computer with $\geq$ 40 GB RAM) is approximately 4.9 hours.

\begin{figure}[htb]
\centering
\includegraphics[width=.7\textwidth]{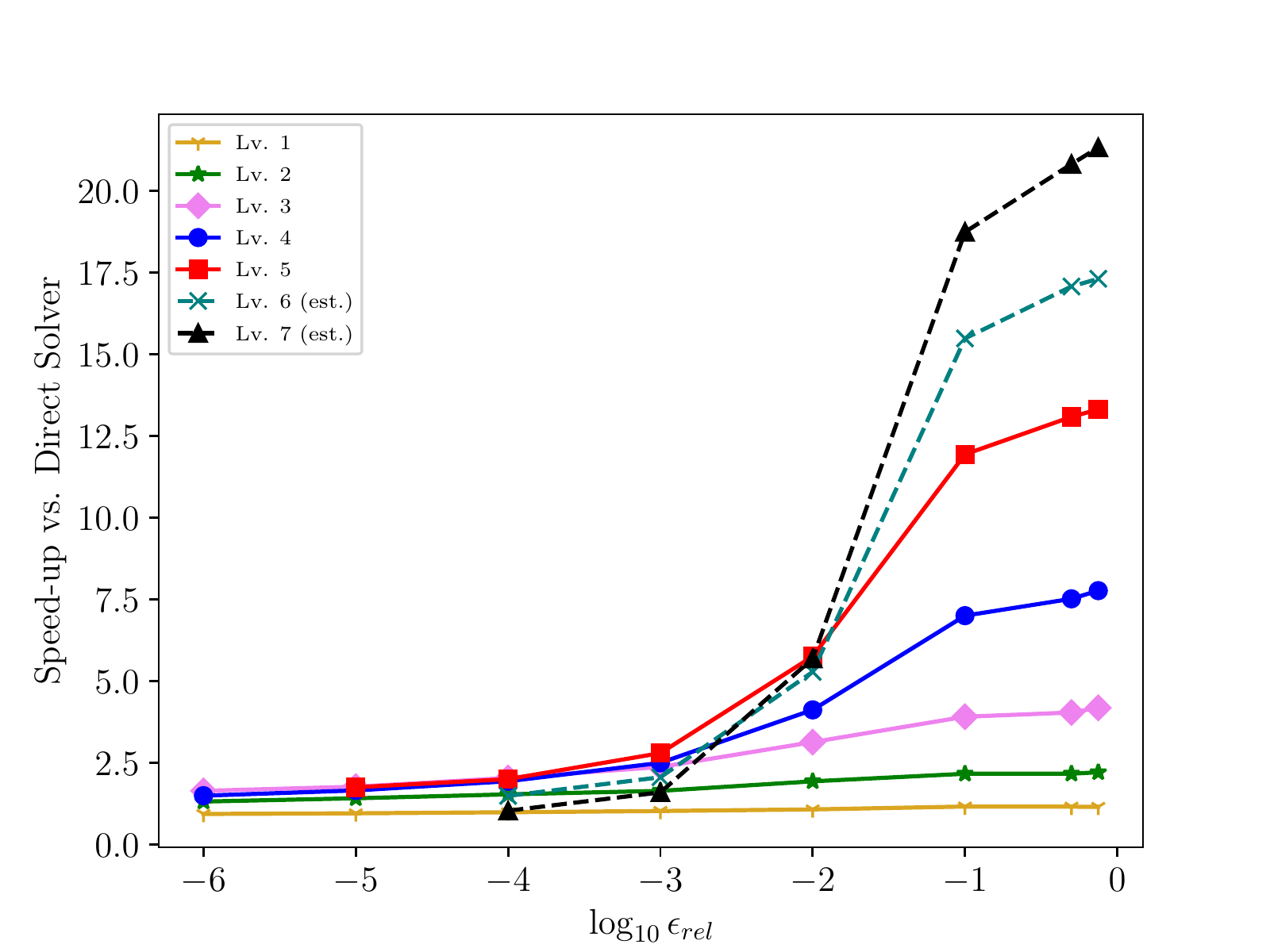}
\caption{Speed-up relative to Direct solve vs.\ $\epsrel$ for different mesh refinement levels. Speed-up for the Level 6 and 7 cases (dashed lines) are based on extrapolated Direct Solve solution times since the compute node had insufficient memory to complete the Direct Solve in those cases.\label{fig:fib_13_spheres_full_piv_ACA_speedup}}
\end{figure}

To conclude this discussion of results, we present in Fig.~\ref{fig:fib_13_spheres_full_piv_ACA_speedup} a plot of the speed-up of the low-rank solver with respect to the direct solver for all meshes and $\epsrel$ values tested. Speed-up here is defined simply as the ratio of the Direct solver time to the low-rank solver time, and in the case of the Level 6 and 7 meshes (dashed lines in the Figure) where the Direct solve time is not available, we use an estimated time based on the curve fit discussed previously.

For the coarsest meshes (Levels 1 and 2), the problem size is too small to realize a large speed-up with respect to the Direct solver; the time spent in ``other'' parts of the code such as I/O and sparse matrix finite element assembly are comparable in both solvers in these cases. On the finest meshes (Levels 6 and 7), the (estimated) speed-up over the Direct solver is approximately \xtimes{17.3} and \xtimes{21.3}, respectively, but we note that the speed-up decreases quickly with decreasing $\epsrel$ (i.e.\ increasing accuracy), so a judicious choice of $\epsrel$ is important for ensuring good performance.

\section{Conclusions\label{sec:conclusions}}

In this work we presented a new approach to the finite element-based computation of nonlinear transient heat transfer with cavity radiation. Conventional finite element formulations for this problem employ dense matrix-based approaches for the cavity terms, which scale poorly as the number of cavity facets becomes large (e.g.\ $> \num{10000}$). To address the poor scaling of the dense matrix formulation, we developed a novel hierarchical low-rank approximation of the cavity terms. As noted in Section~\ref{sec:intro}, this method has similarities to \cite{Potter2022}, but also some important differences, including a low-rank LU-factorization and back/forward substitution framework for efficiently applying cavity terms at many nonlinear iterations and time steps, and an efficient ACA method for constructing low-rank blocks.

Our numerical results demonstrated the accuracy, efficiency, and scalability of the block low-rank framework that we proposed. Using the dense matrix approach as the reference, we were able to obtain highly accurate results in all cases, and with a speed-up of more than \xtimes{20} for the test cases with larger $n_{\rm facets}$. Perhaps most importantly, we demonstrated that for a given $n_{\rm facets}$ the block low-rank approach has much lower memory requirements than the dense approach. This means that, for a given hardware configuration, the block low-rank approach can be applied to much larger models (with higher $n_{\rm facets}$) than is possible with the conventional ``dense matrix'' approach.

 Based on our findings, it is clear that the computational advantage of the block low-rank approach would only increase further for models with larger $n_{\rm facets}$ than those considered here. Hence the methodology proposed in this work is an enabler for solving large-scale cavity radiation problems that arise regularly in industrial and scientific applications, without loss of fidelity in the numerical approximation.



\clearpage
\bibliographystyle{elsarticle-num} 
\bibliography{references}

\begin{thebibliography}{10}
\expandafter\ifx\csname url\endcsname\relax
  \def\url#1{\texttt{#1}}\fi
\expandafter\ifx\csname urlprefix\endcsname\relax\def\urlprefix{URL }\fi
\expandafter\ifx\csname href\endcsname\relax
  \def\href#1#2{#2} \def\path#1{#1}\fi

\bibitem{McQuarrie_1997}
D.~A. McQuarrie, J.~D. Simon, {Physical Chemistry: A molecular approach},
  University Science Books, Sausalito, California, 1997.

\bibitem{Chinoy_1991}
P.~B. Chinoy, D.~A. Kaminski, S.~K. Ghandhi, Effects of thermal radiation on
  momentum, heat, and mass transfer in a horizontal chemical vapor deposition
  reactor, Numerical Heat Transfer, Part A: Applications 19~(1) (1991) 85--100,
  \url{https://doi.org/10.1080/10407789108944839}.

\bibitem{vonZedtwitz_2005}
P.~von Zedtwitz, A.~Steinfeld, {Steam-gasification of coal in a
  fluidized-bed/packed-bed reactor exposed to concentrated thermal radiation
  --- Modeling and experimental validation}, Industrial \& Engineering
  Chemistry Research 44~(11) (2005) 3852--3861,
  \url{https://doi.org/10.1021/ie050138w}.

\bibitem{ZGraggen_2008}
A.~Z'Graggen, A.~Steinfeld, A two-phase reactor model for the
  steam-gasification of carbonaceous materials under concentrated thermal
  radiation, Chemical Engineering and Processing: Process Intensification
  47~(4) (2008) 655--662, \url{https://doi.org/10.1016/j.cep.2006.12.003}.

\bibitem{Barman_2001}
T.~S. Barman, P.~H. Hauschildt, F.~Allard, Irradiated planets, The
  Astrophysical Journal 556~(2) (2001) 885--895,
  \url{https://doi.org/10.1086/321610}.

\bibitem{Delbo_2005}
M.~Delbo, M.~Mueller, J.~P. Emery, B.~Rozitis, M.~T. Capria, Asteroid
  thermophysical modeling, in: P.~Michel, F.~E. Demeo, W.~F. Bottke (Eds.),
  Asteroids, Vol.~IV, The University of Arizona Press, 2005, pp. 107--128,
  \url{https://doi.org/10.2458/azu_uapress_9780816532131-ch006}.

\bibitem{Hayne_2021}
P.~O. Hayne, O.~Aharonson, N.~Sch{\"{o}}rghofer, Micro cold traps on the moon,
  Nature Astronomy 5~(2) (2021) 169--175,
  \url{https://doi.org/10.1038/s41550-020-1198-9}.

\bibitem{Potter2022}
S.~F. Potter, S.~Bertone, N.~Sch{\"{o}}rghofer, E.~Mazarico, Fast hierarchical
  low-rank view factor matrices for thermal irradiance on planetary surfaces,
  ArXiv e-print, \url{https://arxiv.org/abs/2209.07632} (2022).
\newblock \href {http://arxiv.org/abs/2209.07632} {\path{arXiv:2209.07632}}.

\bibitem{Hackbusch2015}
W.~Hackbusch, Hierarchical Matrices: Algorithms and Analysis, Springer, Berlin,
  2015.

\bibitem{Sauter_2011}
S.~A. Sauter, C.~Schwab, Boundary Element Methods, Springer, Berlin, 2011.

\bibitem{Bebendorf2000}
M.~Bebendorf, Approximation of boundary element matrices, Numerische Mathematik
  86~(4) (2000) 565--589, \url{https://doi.org/10.1007/PL00005410}.

\bibitem{Borm_2006}
S.~B{\"{o}}rm, L.~Grasedyck, W.~Hackbusch, {Hierarchical Matrices}, Lecture
  Note No.\ 21, \url{https://tinyurl.com/ycywrvjv} (Jun. 2006).

\bibitem{Franklin_thesis}
A.~M. Franklin, {An implementation of surface-to-surface, blackbody radiation
  heat transfer in a MOOSE application}, Master's thesis, Texas A\&M
  University, \url{https://hdl.handle.net/1969.1/192275} (Aug. 2020).

\bibitem{Holman1986}
J.~P. Holman, Heat Transfer 6th Edition, McGraw-Hill Book Company, Singapore,
  1986.

\bibitem{Chapman_1984}
A.~J. Chapman, Heat Transfer, Fourth Edition, Macmillan, New York, 1984.

\bibitem{Abaqus_2011}
{Abaqus/CAE User's Manual, Version 6.11}, Online (2011).

\bibitem{Dembo_1982}
R.~S. Dembo, S.~C. Eisenstat, T.~Steihaug, {Inexact Newton Methods}, SIAM
  Journal on Numerical Analysis 19~(2) (1982) 400--408,
  \url{https://doi.org/10.1137/0719025}.

\bibitem{Knoll_2004}
D.~A. Knoll, D.~E. Keyes, {Jacobian-free Newton--Krylov methods: A survey of
  approaches and applications}, Journal of Computational Physics 193~(2) (2004)
  357--397, \url{https://doi.org/10.1016/j.jcp.2003.08.010}.

\bibitem{Dunlap_1997}
R.~A. Dunlap, The Golden Ratio and Fibonacci Numbers, World Scientific, New
  Jersey, 1997.

\end{thebibliography}





\end{document}